# TILTED EULER CHARACTERISTIC DENSITIES FOR CENTRAL LIMIT RANDOM FIELDS, WITH APPLICATION TO "BUBBLES"

By N. Chamandy,[1] K. J. Worsley, J. Taylor and F. Gosselin

*McGill University and Google, McGill University and University of Chicago, Stanford University and Université de Montréal*

Local increases in the mean of a random field are detected (conservatively) by thresholding a field of test statistics at a level $u$ chosen to control the tail probability or $p$-value of its maximum. This $p$-value is approximated by the expected Euler characteristic (EC) of the excursion set of the test statistic field above $u$, denoted $\mathbb{E}\varphi(A_u)$. Under isotropy, one can use the expansion $\mathbb{E}\varphi(A_u) = \sum_k \mathcal{V}_k \rho_k(u)$, where $\mathcal{V}_k$ is an intrinsic volume of the parameter space and $\rho_k$ is an EC density of the field. EC densities are available for a number of processes, mainly those constructed from (multivariate) Gaussian fields via smooth functions. Using saddlepoint methods, we derive an expansion for $\rho_k(u)$ for fields which are only approximately Gaussian, but for which higher-order cumulants are available. We focus on linear combinations of $n$ independent non-Gaussian fields, whence a Central Limit theorem is in force. The threshold $u$ is allowed to grow with the sample size $n$, in which case our expression has a smaller relative asymptotic error than the Gaussian EC density. Several illustrative examples including an application to "bubbles" data accompany the theory.

**1. Introduction.** Data which can be modeled as realizations of a smooth random field are increasingly abundant, owing to advances in scanning technology and computer storage. Traditionally, such data have arisen in the fields of oceanography, astronomy and medical imaging (see [40] for an easy introduction). More recently (e.g., [16] or the January 2005 front cover of *Nature* [3]), the "bubbles" paradigm in the cognitive sciences has provided

Received June 2007; revised August 2007.
[1]Supported by a Natural Sciences and Engineering Research Council (NSERC) Canada Graduate Scholarship.
*AMS 2000 subject classifications.* Primary 60G60, 62E20, 62M40; secondary 53A99, 58E05, 60B12, 60F05.
*Key words and phrases.* Saddlepoint approximation, isotropic field, Poisson process, threshold, excursion set.

new and interesting applications. The difficult inference step is most often assessing the significance of the test statistic field arising from such data, as the multiple comparisons problem precludes "pixel-wise" significance testing. Commonly, the supremum of the statistic image is compared to an estimate of its upper $p$th quantile under the global null hypothesis of a zero mean function (equivalently, the statistic field is *thresholded* at this quantile to identify regions of activation). A leading approach to approximating these quantiles uses the expected Euler characteristic (EC) of the excursion of the field above a threshold ([2, 38]), as described below. However, a common thread in these analyses is the assumption that the random field data are Gaussian; in some applications this assumption is known to be flawed.

Practitioners of voxel-based morphometry in neuroimaging have highlighted the problem [27, 37], but have stopped short of proposing a theoretical solution. In the context of galaxy density random fields, Matsubara [21, 22] obtained approximations to the expected EC in three dimensions using Edgeworth series (though he never made the jump to a saddlepoint expansion). Catelan and coauthors [8] had previously used a similar approach (also in an astrophysics context) to get at the mean cluster size of excursions for non-Gaussian fields. Rabinowitz and Siegmund [25] used saddlepoint methods to approximate the tail probability for the supremum of a particular non-Gaussian field, a smoothed Poisson point process. Although they supplied a first-order version of the expected EC and published heuristic arguments, their work originated some of the central lines of reasoning of the current article (indeed, since the paper is cited often we henceforth refer to it as RS). We shall formulate a general approximation to the expected EC for asymptotically Gaussian fields with known (decaying) cumulants, and compare its relative error to that under a full Gaussian assumption. In particular, our approximation is found to outperform the Gaussian expected EC in a particular "large threshold" asymptotic setting. Throughout, we underscore the connection between our results and those of RS.

For concreteness, we consider processes of the form

$$(1.1) \qquad Z(t) = \frac{1}{\sqrt{n}} \sum_{i=1}^{n} W_i(t),$$

where the $W_i$'s are i.i.d. smooth random fields on $T \subset \mathbb{R}^N$ with zero mean, constant variance $\sigma^2$ and finite higher cumulants. Such a $Z(t)$ is a special case of what we call a Central Limit random field (CLRF), and usually should be thought of as a test statistic derived from data. We call $A_u = A_u(Z,T) = \{t \in T; Z(t) \geq u\}$ the *excursion set* of $Z$ above $u$. We denote the EC, an integer-valued set functional, by $\varphi(\cdot)$. For a working definition of the EC, see [34], and for a more rigorous account the reader is referred to Part II of Adler and Taylor [2]. As the latter volume is the second major source of methodology used in this article, it will be referenced as AT. Very loosely (ignoring what happens on the boundary), in 2D the EC counts #"blobs" − #"holes" in a set, while in 3D it counts #"blobs" − #"tunnels" + #"hollows" [40].

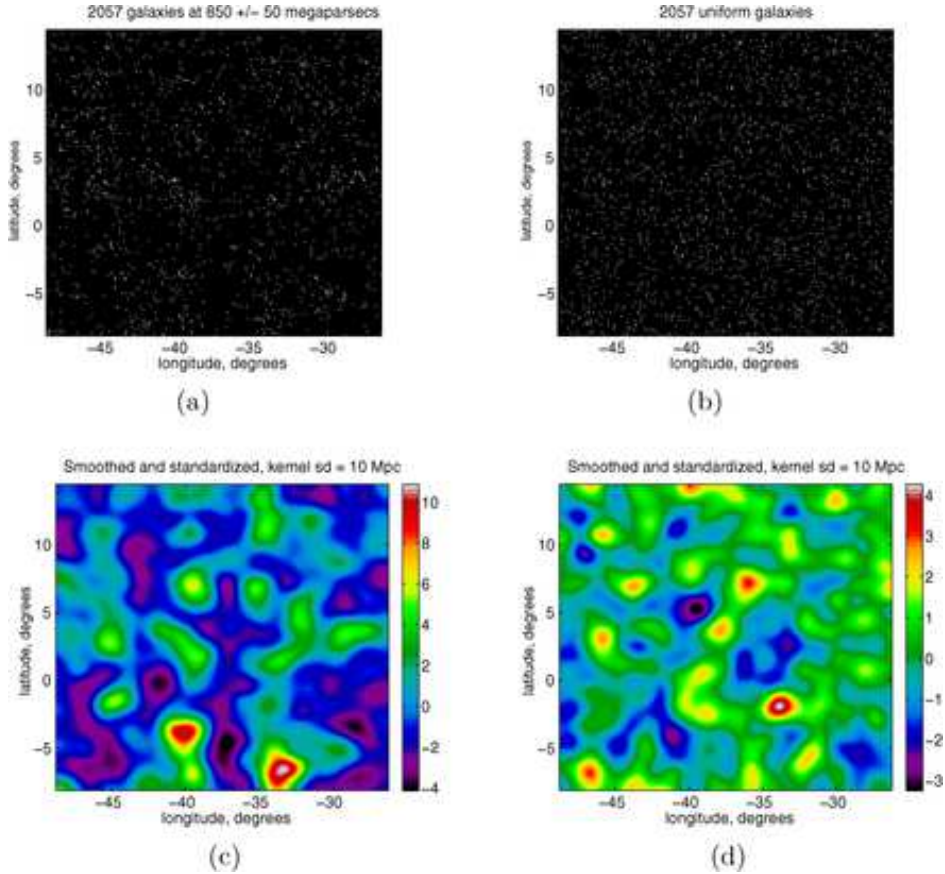

FIG. 1. *Galaxies are not uniformly distributed in the universe, but tend to congregate in clusters, strings or even sheets. Cosmologists have used the EC of regions of high galaxy density to characterize this large-scale structure ([18] used data from the first release of the Sloan Digital Sky Survey SDSS). In 2D the EC counts #"blobs"—#"holes" in a set. We examine luminous red galaxies (LRG's—the most distant objects in the SDSS) from the latest (5th) data release [24]. In* (a) *we show LRG's from a* 22.5° *square patch at* 850 ± 50 *megaparsecs (Mpc) with roughly the same galaxy density as in [18] (our sample is part of a sphere, rather than cone-shaped, so we avoid the luminosity correction). Panel* (b) *displays the same number of "galaxies" simulated uniformly according to a Poisson process. In* (c)–(d) *the LRG's are smoothed with a Gaussian kernel of s.d. 10 Mpc, and standardized to have zero mean and unit variance under the Poisson model.*

In some applications the expected EC of $A_u$ as a function of $u$ is itself of central interest. For example, cosmologists compare the observed EC (or $1 - $ EC, sometimes called the *genus*) from galaxy survey data with the expected EC from various hypothesized models of the large-scale structure of the universe (see [18] and the illustration in Figure 1). In other applications, including neuroimaging and the bubbles problem treated in Section 8, interest centers around using this quantity as an approximation to the $p$-value:

$$(1.2) \qquad \mathbb{E}\varphi(A_u) \approx \mathbb{P}\left\{\sup_{t \in T} Z(t) \geq u\right\}$$

for large $u$. The assumption that (1.2)—the so-called "EC heuristic"—is adequate for statistical inference is no longer a mere heuristic for Gaussian fields; see [33] or AT, Chapter 14. See also Section 4.4 of the present paper.

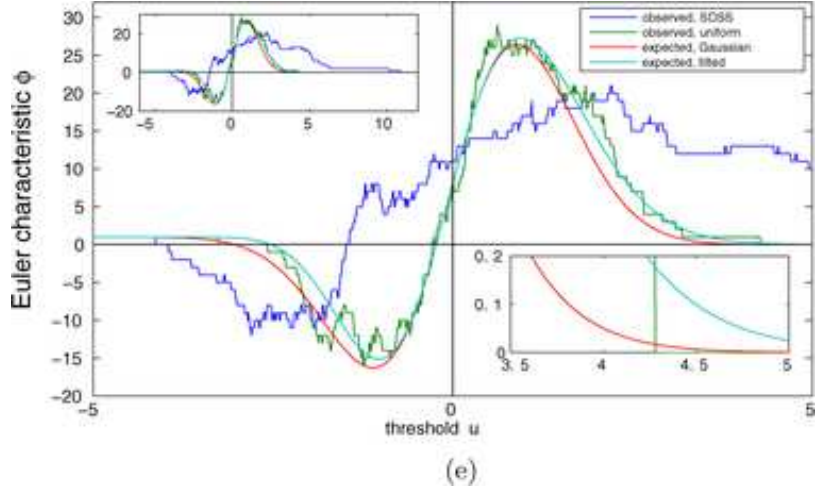

FIG. 1. *(Continued.) In* (e) *we compare the observed EC of each image at a continuum of thresholds with the expected EC under the Poisson model computed using a Gaussian approximation and the "tilted" correction derived in Sections 4–5. The observed EC for the SDSS image deviates from the others, reinforcing the nonuniformity of LRG's. The Gaussian and tilted expected EC's differ primarily in the tails; in the range where the expected EC approximates a p-value (bottom-right zoom) this correction is dramatic. In particular, the tilted curve majorizes the Gaussian one since the LRG distribution has positive cumulants. The observed EC for the (asymptotically Gaussian) Poisson data conforms well with both curves. Note that all curves approach 1 as $u \to -\infty$ and 0 as $u \to \infty$ (top-left zoom).*

To motivate the study of what are known as *EC densities*, we have the following well-established result.

THEOREM 1 ([39], AT).    *Let $Z(t)$ be any isotropic random field on $T \subset \mathbb{R}^N$, a locally convex finite union of convex bodies. Then the expected EC of* *its excursion set at the level $u \in \mathbb{R}$ can be written*

$$\text{(1.3)} \qquad \mathbb{E}\varphi(A_u) = \sum_{k=0}^{N} \mathcal{V}_k(T)\rho_k(u),$$

*where $\mathcal{V}_k(T)$ is the kth intrinsic volume of the search region and $\rho_k(u)$ is a function called the kth EC density of the field.*

PROOF. The expected EC $\mathbb{E}\varphi(A_u(Z,\cdot))$ is a real-valued additive functional acting on $T$; $\mathbb{E}\varphi(A_u(Z, T \cup V)) = \mathbb{E}\varphi(A_u(Z,T) \cup A_u(Z,V)) = \mathbb{E}\varphi(A_u(Z,T)) + \mathbb{E}\varphi(A_u(Z,V)) - \mathbb{E}\varphi(A_u(Z,T) \cap A_u(Z,V)) = \mathbb{E}\varphi(A_u(Z,T)) + \mathbb{E}\varphi(A_u(Z,V)) - \mathbb{E}\varphi(A_u(Z,T \cap V))$. If the process is isotropic, then $\mathbb{E}\varphi(A_u(Z,\cdot))$ is also translation- and rotation-invariant. Finally, the local convexity of $T$ guarantees that the functional is continuous (e.g., in the Hausdorff metric). The result then follows from Hadwiger's seminal theorem [17, 28], which states that any such functional is in the linear span of the intrinsic volumes, and hence admits an expansion like (1.3). □

REMARK. If the field is nonisotropic but derived from a (multivariate) Gaussian field by a smooth function, then on sufficiently regular parameter spaces (1.3) holds with $\mathcal{V}_k(T)$ replaced by the $k$th Lipschitz–Killing curvature $\mathcal{L}_k(T)$. The latter differential geometric functionals depend both on the space $T$ and on the covariance structure of the field itself; see AT Chapter 10.

1.1. *Intrinsic volumes.* The intrinsic volumes of a set $T$ are a generalization of its volume to lower dimensional measures. There are several ways of describing these functionals; here we give an implicit definition, which is valid for any locally convex set $T$ and sufficiently small radius $r$. Let $|\cdot|$ denote the Lebesgue measure, $B(0,r)$ be the ball of radius $r$ centred at 0, and $T \oplus B(0,r) = \{x + y : x \in T, y \in B(0,r)\}$ be the tube of radius $r$ around $T$. Write $\omega_k = |B_k(0,1)| = \pi^{k/2}/\Gamma(k/2+1)$ for the Lebesgue measure of the unit ball in $\mathbb{R}^k$. Then

$$\text{(1.4)} \qquad |T \oplus B(0,r)| = \sum_{k=0}^{N} \omega_{N-k} r^{N-k} \mathcal{V}_k(T).$$

Thus, the intrinsic volumes are related to the coefficients in the Steiner–Weyl volume of tubes expansion. For example, in $\mathbb{R}^3$ we have $\mathcal{V}_0(T) = \varphi(T)$, $\mathcal{V}_1(T) = 2 \times \text{caliper diameter}(T)$, $\mathcal{V}_2(T) = 1/2 \times \text{surface area(T)}$, and $\mathcal{V}_3(T) = \text{volume}(T)$. Remarkably, EC densities for all "Gaussian-related" processes were shown to have a similar characterization in [31].

1.2. *Computing EC densities.* For a matrix $A$ or vector $a$ we introduce the notation $A|_k = (A_{ij})_{1 \le i,j \le k}$ and $a|_k = (a_i)_{1 \le i \le k}$. When $k = 0$ this denotes, by convention, an "empty" matrix or vector—that is, a term or factor which disappears from the expression. We use $I_p$ (or just $I$ when the dimension is unambiguous) to denote the $p \times p$ identity matrix. When applied to a random field, "˙" represents spatial differentiation with respect to $t$, and we shall refer to the individual components of $\dot{Z}$ and $\ddot{Z}$, respectively, by $Z_i$ and $Z_{ij}$. There are several useful formulations for computing the EC densities of stationary fields [39]. In the derivation of Theorem 4 we use

$$(1.5) \quad \rho_k(u) = \mathbb{E}\{Z_k \mathbb{1}_{[Z_k \ge 0]} \det(-\ddot{Z}|_{k-1})|Z = u, \dot{Z}|_{k-1} = 0\} f_{0,k-1}(u, 0)$$

for $k = 1, \ldots, N$, where $f_{0,k-1}(\cdot)$ denotes the joint density of $(Z, \dot{Z}|'_{k-1})$. Indeed, this was the original method used by Adler to compute the expected EC, and corresponds to using the $N$th coordinate function $(t_1, t_2, \ldots, t_N) \mapsto t_N$ (instead of the field itself) as the function in Morse's Theorem (see [1], Chapter 4). The zeroth EC density is just the univariate tail probability, $\rho_0(u) \stackrel{\Delta}{=} \mathbb{P}\{Z(t) \ge u\}$. It has a clear interpretation as the lowest order approximation to the $p$-value (1.2).

1.3. *A glimpse ahead.* EC densities are available in diverse parameter spaces for a number of processes [6, 7, 30, 32, 35, 38], mainly those derived from (multivariate) Gaussian fields via smooth functions. Let $H_m(x) = m! \sum_{j=0}^{\lfloor m/2 \rfloor} \frac{(-1)^j x^{m-2j}}{2^j j!(m-2j)!}$ be the $m$th Hermite polynomial. For a mean zero Gaussian field with variance $\sigma^2$, the densities take the form (cf. AT, Chapter 11)

$$(1.6) \qquad \rho_k^\gamma(u) \stackrel{\Delta}{=} (2\pi)^{-(k+1)/2} \sigma^{-k} \exp\left\{-\frac{u^2}{2\sigma^2}\right\} H_{k-1}\left(\frac{u}{\sigma}\right).$$

(Note that $\gamma$ is used to denote the Gaussian measure, following the convention of AT.) Our main result (Theorem 4) will be to show that for a stationary CLRF like (1.1), and assuming that the threshold $u$ is of the order $n^{1/2-\alpha}$ for some $1/4 < \alpha \le 1/2$, modifying only the exponential portion of (1.6) leads to an approximation to $\rho_k(u)$ whose relative error behaves like $n^{1/2-2\alpha}$. Adding factor related to the mixed skewness terms $\mathbb{E}[Z^2 Z_i]$ can yield an additional improvement. Moreover, we shall show that this error is asymptotically strictly better than that of (1.6), exponentially so over the range $1/4 < \alpha \le 1/3$.

The remainder of the paper is organized as follows. In Section 2 we use saddlepoint methods to approximate the joint density of $(Z, \dot{Z}', \ddot{Z}')$ (suitably transformed), leading to a mixed tilted-direct Edgeworth expansion. In Section 3 we consider the growth rate of $u$ and compute certain moments of the tilted distribution which shall be used in subsequent derivations. We derive the approximation to $\rho_k(u)$, by methods reminiscent of [1] and [38], in Section 4. Section 5 offers an analysis of the relative error of the tilted densities compared to that of a standard Gaussian approximation, and we comment on the contribution of RS in Section 6. In Section 7 we illustrate the method by treating a field—the $\chi_n^2$ process—whose densities are known exactly, and Section 8 applies the theory to the test statistic from a bubbles experiment. We also present a simulation under the latter framework which examines the accuracy of our corrected version of (1.3) as a $p$-value approximation. In the interest of brevity, we have omitted several proofs and technical developments, all of which can be found in an online version of the present article (URL), or in the doctoral thesis [9].

**2. A mixed expansion.** In order to avoid having EC densities which depend on $t$ in Theorem 4 below, we shall assume that $Z$ is strictly stationary. For practical purposes [i.e., to apply Theorem 1 in the computation of $\mathbb{E}\varphi(A_u)$], the stronger requirement that $Z$ be strictly isotropic shall be in force. In that case, purely for convenience, $Z$ will be assumed to have unit second spectral moment, namely $\mathbb{E}[Z_i Z_j] = \delta_{ij}$. The latter assumption can easily be relaxed to allow $\mathsf{Var}[\dot{Z}(t)] \equiv \lambda I$, which merely modifies (1.3) by multiplying the $k$th term by $\lambda^{k/2}$. Strict isotropy can, on very simple parameter spaces, be weakened to stationarity [with a similar correction to (1.3)]; see Section 4.3 or AT, Chapter 11. If $\xi \geq 0$, the notation $x_n = O(n^{-\xi})$ will indicate that for some $0 < C < \infty$ and all sufficiently large $n$, $|x_n| \leq C n^{-\xi}$. A similar notation is used for matrices, with the understanding that the relation holds componentwise. The notation $a_n \sim b_n$ signifies that $a_n/b_n \to 1$ as $n \to \infty$.

Fix $k \geq 2$ and put

$$(2.1) \qquad M = M_k \triangleq \sigma \ddot{Z}|_{k-1} + \frac{Z}{\sigma} I_{k-1}.$$

The most important property of this shifted $(k-1) \times (k-1)$ Hessian is that it is uncorrelated with the field: $\mathbb{E}\{M_{ij} Z\} \equiv 0$. Let $\lambda_{ijlm} \triangleq \mathbb{E}\{Z_{ij} Z_{lm}\}$ [a symmetric, $O(1)$ function] denote a fourth-order spectral moment of the process. It is easily shown that

$$\lambda_{ijlm}^M \triangleq \mathbb{E}\{M_{ij} M_{lm}\} = \sigma^2 \lambda_{ijlm} - \delta_{ij}\delta_{lm}.$$

We shall write $M$ to denote both the matrix and the corresponding random vector $\mathrm{vec}(M)$, with the context eliminating confusion (similarly for $\ddot{Z}$). In view of this, let $\mathbf{Z} = \mathbf{Z}_k^n \triangleq (Z, \dot{Z}|_k', M')'$ have density $f(\mathbf{z})$, to which we shall apply the asymptotic expansion. Define the exponential family

$$(2.2) \qquad f_\theta(\mathbf{z}) = f_\theta(z, \dot{z}, m) = \exp\{\theta \mathbf{z} - K(\theta)\} f(z, \dot{z}, m),$$

where $\theta = (\theta_0, \ldots, \theta_{d_k-1}) \in \mathbb{R}^{d_k} = k(k+1)/2 + 1$. The function $K(\cdot)$ is of course the cumulant generating function (cgf) of $\mathbf{Z}$,

$$K(\theta) = \sum_{\nu=1}^{\infty} \sum_{\substack{\gamma \in \mathbb{N}^{d_k} \\ |\gamma|=\nu}} \frac{1}{\gamma!} k_\nu([\mathbf{Z}]^\gamma) \theta^\gamma, \tag{2.3}$$

where if $\gamma \in \mathbb{N}^d$, we define $[Y]^\gamma = [\overbrace{Y_1, \ldots, Y_1}^{\gamma_1}, \ldots, \overbrace{Y_d, \ldots, Y_d}^{\gamma_d}]$, $(\theta_1, \ldots, \theta_d)^\gamma = \theta_1^{\gamma_1} \cdots \theta_d^{\gamma_d}$, $|\gamma| = \sum_{i=1}^{d} \gamma_i$, and $\gamma! = \prod_{i=1}^{d} \gamma_i!$. We have written $k_\nu$ for a $\nu$th order cumulant.

Tilting in the first argument and inverting (2.2) yields

$$f(z, \dot{z}, m) = \exp\{K(\theta_0, 0') - \theta_0 z\} f_{(\theta_0, 0')}(z, \dot{z}, m). \tag{2.4}$$

The relationship (2.4) holds for all values of $\theta_0$, in particular, that value $\theta_0 = \hat{\theta}_u$, which makes the Edgeworth expansion to $f_{(\theta_0, 0')}(\mathbf{z})$ most accurate at the point $\mathbf{z} = (u, 0', 0')'$—namely, that for which

$$u = \mathbb{E}_{(\hat{\theta}_u, 0')} Z = \left.\frac{\partial K(\theta_0, 0')}{\partial \theta_0}\right|_{\theta_0 = \hat{\theta}_u} \tag{2.5}$$

holds. It is trivial to see that this requirement is satisfied if and only if $\hat{\theta}_u$ is the formal maximum likelihood estimate of $\theta_0$ under the model $f_{(\theta_0, 0')}(\mathbf{z})$ when $Z = u$ is observed.

Finally, we let $\mu_/$ and $\Sigma_/$ denote the mean and covariance matrix under the density $f_{(\hat{\theta}_u, 0')}(\mathbf{z})$. Then (cf. [4], page 187) we have

$$\mu_/ = \begin{pmatrix} u \\ \nabla_{2\cdots d_k} K(\theta)|_{\theta=(\hat{\theta}_u, 0')} \end{pmatrix} \quad \text{and} \quad \Sigma_/ = \nabla\nabla' K(\theta)|_{\theta=(\hat{\theta}_u, 0')},$$

where $\nabla_{2\cdots d_k}$ has been used to denote $(\partial/\partial\theta_2, \ldots, \partial/\partial\theta_{d_k})'$. The approximation to be exploited will eventually take the form

$$\begin{aligned}
f(z, \dot{z}, m) &= \exp\{K(\hat{\theta}_u, 0') - \hat{\theta}_u z\} \phi_{d_k}((z, \dot{z}', m')'; \mu_/, \Sigma_/) \\
&\quad \times \{1 + Q_{k,3}((z, \dot{z}', m')' - \mu_/)\} + e_{k,n}, \\
&\triangleq \hat{f}_k(z, \dot{z}, m)\{1 + Q_{k,3}((z, \dot{z}', m')' - \mu_/)\} + e_{k,n},
\end{aligned} \tag{2.6}$$

where $Q_{k,3}$ is a degree 3 polynomial containing mixed cumulants of $(Z, \dot{Z}|_k, M)$ of order 2 and 3, while $e_{k,n}$ denotes terms of smaller asymptotic order. It turns out that $Q_{k,3}$ has coefficients (including constant term) which range up to $O(n^{-1/2})$. Expansion (2.6) is called "mixed tilted-direct" because exponential tilting is applied to the first argument but not to the others.

**3. Preliminary calculations.** We argued heuristically that the expected EC will provide a good approximation to the $p$-value (1.2) for large thresholds $u$; we therefore allow the possibility that $u = u_n \to \infty$ as $n \to \infty$. Let $u = O(n^{1/2-\alpha})$ for some $0 \leq \alpha \leq 1/2$. From (2.5) we have

$$
\begin{aligned}
u_n &= K'_Z(\hat{\theta}_u) \\
&= \sigma^2 \hat{\theta}_u + \frac{1}{2} k_3(Z) \hat{\theta}_u^2 + \frac{1}{3!} k_4(Z) \hat{\theta}_u^3 + \cdots,
\end{aligned}
\tag{3.1}
$$

where $K_Z(\cdot)$ denotes the marginal cgf of $Z$. When $\theta$ is a univariate argument, we shall use "′" or a superscript [e.g., "$(3)$"] to denote (multiple) differentiation with respect to $\theta$. Under some additional assumptions, (3.1) is a proper asymptotic expansion, in the sense that the error in the remainder is of smaller order than the previous term. We shall use the following lemma repeatedly, and often implicitly, in proving our main results.

LEMMA 2. *Let $Z = n^{-1/2} \sum_i W_i$ with $W_i$ centred, i.i.d. random variables with variance $\sigma^2$ such that $\mathbb{E}|W_i|^{s+1} < \infty$ for some $s \geq 1$. Suppose that $K_W^{(s+1)}$ is continuous on a closed interval $[-\epsilon_s, \epsilon_s]$ with $\epsilon_s > 0$. Let $|\theta| \leq cn^{1/2-\beta}$ for some finite $c$ and $\beta > 0$. Then there exists a constant $c(s)$ depending only on $s$ such that, for all large $n$,*

$$
\left| K'_Z(\theta) - \sigma^2 \theta - k_3(Z) \frac{\theta^2}{2} - \cdots - k_s(Z) \frac{\theta^{s-1}}{(s-1)!} \right| \leq c(s) n^{1/2-s\beta}.
\tag{3.2}
$$

*In particular, if $u \sim cn^{1/2-\alpha}$ for some $0 < \alpha \leq 1/2$ and $c > 0$, with $Z$ as described above and $s = 2$, then*

$$
\hat{\theta}_u = \hat{\theta}_{u,n} = \frac{u}{\sigma^2}[1 + O(n^{-\alpha})]
\tag{3.3}
$$

*as $n \to \infty$. If $\alpha = 1/2$, then $u$ is (asymptotically) constant and $\hat{\theta}_u = u/\sigma^2 + O(n^{-1/2})$.*

3.1. *The tilted moments.* Let us now be more explicit about the "tilted" parameters $\mu_/$ and $\Sigma_/$, computed by differentiating (2.3) with respect to the appropriate argument(s) and then setting $\theta = (\hat{\theta}_u, 0')$. To simplify notation, we partition as follows:

$$
\mu_/ = (u, \dot{\mu}', \ddot{\mu}')'; \qquad \Sigma_/ = \begin{pmatrix} \tau^2 & \dot{\sigma}' & \ddot{\sigma}' \\ \dot{\sigma} & \dot{\Sigma} & \ddot{\Sigma}' \\ \ddot{\sigma} & \ddot{\Sigma} & \ddddot{\Sigma} \end{pmatrix}.
$$

We alert the reader that the dot notation in the above parameters does not indicate differentiation, but was designed to emphasize that portion of the process $Y(t)$ with which a particular mean or covariance is associated. In $\Sigma_/$, the diagonal entries in the partition are just the tilted variances, respectively, of $Z$, $\dot{Z}|_k$ and $M$. The off-diagonal entries are the corresponding covariances, and so $\dot{\sigma}$ is $k \times 1$, $\ddot{\sigma}$ is $k(k-1)/2 \times 1$, and $\dddot{\Sigma}$ is $k(k-1)/2 \times k$. We have, for example,

$$\dot{\mu}_i = \sum_{\nu=3}^{\infty} \frac{1}{(\nu-1)!} k_\nu(Z, \ldots, Z, Z_i) \hat{\theta}_u^{\nu-1},$$

(3.4)

$$\dddot{\Sigma}_{i,jl} = \sum_{\nu=3}^{\infty} \frac{1}{(\nu-2)!} k_\nu(Z, \ldots, Z, Z_i, M_{jl}) \hat{\theta}_u^{\nu-2},$$

and many more such expressions. When $\alpha > 1/4$ it is easily deduced that

$$\dot{\mu}_i = O(n^{1/2-2\alpha}), \qquad \ddot{\mu}_{ij} = O(n^{1/2-2\alpha});$$

(3.5) $\quad \tau_u^2 \triangleq \tau^2 = \sigma^2 + O(n^{-\alpha}), \qquad \dot{\sigma}_i = O(n^{-\alpha}), \qquad \ddot{\sigma}_{ij} = O(n^{-\alpha});$

$$\dot{\Sigma}_{ij} = \delta_{ij} + O(n^{-\alpha}), \qquad \ddot{\Sigma}_{ij,lm} = \lambda^M_{ijlm} + O(n^{-\alpha}),$$

$$\dddot{\Sigma}_{i,jl} = O(n^{-\alpha}).$$

3.2. *Some tilted conditional moments.* We shall also need the following conditional moments (assuming $\alpha > 1/4$) under the tilted Normal density $\phi_k(\cdot) \triangleq \phi_{d_k}(\cdot; \mu_/, \Sigma_/)$, the proofs of which are straightforward:

(3.6) $\qquad \nu \triangleq \mathbb{E}_{\phi_k}\{Z_k | Z = u, \dot{Z}|_{k-1} = 0\} = O(n^{1/2-2\alpha});$

(3.7) $\qquad \eta^2 \triangleq \mathsf{Var}_{\phi_k}\{Z_k | Z = u, \dot{Z}|_{k-1} = 0\} = 1 + O(n^{-\alpha});$

(3.8)
$$\mathbb{E}_{\phi_k}\{M_{ij} | Z = u, \dot{Z}|_{k-1} = 0, Z_k = x\}$$
$$= \ddot{\mu}_{ij} + c_{ij}x + O(n^{1/2-3\alpha});$$

(3.9)
$$\mathbb{E}_{\phi_k}\{M_{ij}M_{lm} | Z = u, \dot{Z}|_{k-1} = 0, Z_k = x\}$$
$$= \lambda^M_{ijlm} + c^0_{ijlm} + c^1_{ijlm}x + c^2_{ijlm}x^2;$$

where $c_{ij} = O(n^{-\alpha})$ and

(3.10)
$$c^0_{ijlm} = O(n^{\max(-\alpha, 1-4\alpha)}),$$
$$c^1_{ijlm} = O(n^{1/2-3\alpha}),$$
$$c^2_{ijlm} = O(n^{-2\alpha}).$$

The punchline to (3.8) and (3.9) is that under the Gaussian approximation to the conjugate density $f_{(\hat{\theta}_u, 0')}$, and conditional on the event $\{Z = u, \dot{Z}|_{k-1} =$ $0, Z_k = x\}$, we have $M \stackrel{\mathscr{D}}{=} \Delta + \Xi_x$, where $\Delta$ is a centred Gaussian matrix with covariances satisfying (4.1) below, and $\Xi_x$ is a small, nonrandom (given $x$) perturbation.

**4. Euler characteristic densities.** We shall now proceed to use the characterization (1.5) in order to derive our corrected EC density.

4.1. *The case $\alpha > \frac{1}{4}$.* The following analogue of AT's Lemma 11.5.2 can be proved in a similar fashion if one is mindful of the various error terms. In what follows, $|\cdot|$ is the determinant. We also add a subscript $n$ to all EC density expressions, since we are interested in their asymptotic behavior.

LEMMA 3. *Let $\frac{1}{4} < \alpha \leq \frac{1}{2}$, let $0 \leq x$ be a $O(1)$ nonrandom scalar and $\varepsilon$ a $O(1)$ symmetric function from $\{1,\ldots,k\}^4$ into $\mathbb{R}$. Suppose that $\Delta$ is a centered, (jointly) Gaussian $k \times k$ matrix whose covariances satisfy*

$$(4.1) \qquad \mathbb{E}(\Delta_{ij}\Delta_{lm}) = \varepsilon(i,j,l,m) - \delta_{ij}\delta_{lm} + s_{ijlm}(x),$$

*where $s_{ijlm}(x)$ is a quadratic polynomial with coefficients as in (3.10). Let $b = O(n^{1/2-\alpha})$ be a nonrandom scalar, and let $\Xi^0 = O(n^{1/2-2\alpha})$, $\Xi^1 = O(n^{-\alpha})$, and $\Xi_x = \Xi^0 + x\Xi^1$ be nonrandom matrices. Then*

$$(4.2) \qquad \mathbb{E}|\Delta - bI + \Xi_x| = (-1)^k H_k(b)\{1 + r_k(x)\},$$

*where $r_k$ is a degree $k$ polynomial whose largest (constant) coefficient is $a_0 = O(n^{-\alpha})$ as $n \to \infty$.*

THEOREM 4. *Let $Z$ be a strictly stationary CLRF as in (1.1), with covariance function belonging to $C^4(\mathbb{R}^N)$ and $\mathbb{E}[\dot{Z}\dot{Z}'] = I$. Let $u \sim cn^{1/2-\alpha}$ as $n \to \infty$ for some constants $0 < c < \infty$ and $\frac{1}{4} < \alpha \leq \frac{1}{2}$, and suppose that $\hat{\theta}_u$ is a solution to*

$$(4.3) \qquad K'_Z(\theta) = u,$$

*where $K_Z(\cdot)$ is the cgf of $Z(t)$. Let $\tau_u^2 = K''_Z(\hat{\theta}_u)$ be the tilted variance of $Z(t)$ at $\hat{\theta}_u$, and $I_u(\theta) = u\theta - K_Z(\theta)$. Then subject to regularity conditions, the EC densities of $Z$ for $k \geq 1$ are given by*

$$(4.4) \qquad \rho_{k,n}(u) = \widehat{\rho}_{k,n}(u) \times \{1 + \epsilon_n\},$$

*where*

$$(4.5) \qquad \widehat{\rho}_{k,n}(u) = (2\pi)^{-(k+1)/2} \tau_u^{-k} \exp\{-I_u(\hat{\theta}_u)\} H_{k-1}(\tau_u\hat{\theta}_u)$$

*and as $n \to \infty$ the relative error satisfies $\epsilon_n = O(n^{1/2-2\alpha})$.*

*Moreover, an approximation with relative error $O(n^{\max(-\alpha, 1-4\alpha)})$ is given by*

$$\widehat{\rho}_{k,n}(u) \times \left(1 + \mathbb{E}[Z^2 Z_k]\sqrt{\frac{\pi}{2}}\frac{\hat{\theta}_u^2}{2}\right), \tag{4.6}$$

*and (if $k \geq 2$) one with relative error $O(n^{\max\{-\alpha, 3/2 - 6\alpha\}})$ by*

$$\widehat{\rho}_{k,n}(u) \times \left(1 + \mathbb{E}[Z^2 Z_k]\sqrt{\frac{\pi}{2}}\frac{\hat{\theta}_u^2}{2}\right)\exp\left\{-\frac{\hat{\theta}_u^4}{8}\sum_{i=1}^{k-1}\mathbb{E}[Z^2 Z_i]^2\right\}. \tag{4.7}$$

REMARK. The "regularity conditions" mentioned above are just those on the random field $(Z, \dot{Z}|_k', \ddot{Z}|_{k-1}')$ required to apply the Expectation Meta-Theorem (Theorem 11.1.1 in AT), those on the density $f(z, \dot{z}, m)$ assuring the correctness of the Edgeworth expansion, plus a mild assumption about the moments of $(W, \dot{W}, \ddot{W})$. These conditions are academic in that they can rarely be checked in practice, but are considerably weaker than the Gaussian assumption, which is currently the only alternative. See the Appendix for further details.

PROOF OF THEOREM 4 (LEADING TERM). For extra details the reader is referred to the online version of this article, where we rigorously bound the error terms and consider the contribution to the integral (1.5) of the polynomial $Q_{k,3}(\cdot)$. It is seen there that the additional terms are of relative size at most $O(n^{-\alpha})$, and are therefore negligible. The crux of the derivation, however, is in working out the leading term, which we do now. Plugging the saddlepoint approximation into (1.5) yields

$$\rho_{k,n}(u) \doteq \int_{x=0}^{\infty}\int_{m \in \mathbb{R}^{D_{k-1}}} x \det\left(-\frac{m - u/\sigma I}{\sigma}\right)\widehat{f}_k(u, 0', x, m)\, dm\, dx \tag{4.8}$$

$$= \exp\{-I_u(\hat{\theta}_u)\}\left(\frac{-1}{\sigma}\right)^{k-1}$$

$$\times \int_{x=0}^{\infty} x\phi_k(x|u, 0')$$

$$\times \int_{\mathbb{R}^{D_{k-1}}} \det(m - u/\sigma I)\phi_k(m|u, 0', x)\, dm\, dx \tag{4.9}$$

$$\times \phi_{0,k-1}(u, 0'),$$

where $D_j = j(j+1)/2$ and the various $\phi$'s denote the obvious Gaussian densities derived from $\phi_{d_k}(\cdot; \mu_/, \Sigma_/)$ by conditioning or marginalizing. Evaluating the inner expectation via (3.8)–(3.9) and Lemma 3, after noting that

$\tau_u \hat{\theta}_u - u/\sigma = O(n^{-\alpha})$ has smaller order than $\Xi^0$, gives

$$\exp\{-I_u(\hat{\theta}_u)\}\sigma^{-(k-1)}H_{k-1}(\tau_u\hat{\theta}_u)$$
(4.10)
$$\times \mathbb{E}\{Y(1+r_{k-1}(Y))\mathbb{1}_{[Y>0]}\}\phi_{0,k-1}(u,0'),$$

where $Y \sim \mathcal{N}(\nu,\eta^2)$ as in (3.6) and (3.7). Recall that $\sigma = \tau_u(1+O(n^{-\alpha}))$, and consider the last factor in (4.10). We have

(4.11) $\quad \phi_{0,k-1}(u,0') = (2\pi)^{-k/2}\det(V_k)^{-1/2}\exp\{-\tfrac{1}{2}\dot{\mu}|'_{k-1}\dot{\Sigma}|^{-1}_{k-1}\dot{\mu}|_{k-1}\},$

where

(4.12) $$V_k \triangleq \begin{pmatrix} \tau_u^2 & \dot{\sigma}|'_{k-1} \\ \dot{\sigma}|_{k-1} & \dot{\Sigma}|_{k-1} \end{pmatrix}.$$

So $\det(V_k) = \tau_u^2(1+O(n^{-\alpha}))$, and the exponential portion of (4.11) behaves like

$$\exp\{O(n^{1/2-2\alpha})'[I+O(n^{-\alpha})]O(n^{1/2-2\alpha})\} = \exp\{O(n^{1-4\alpha})\}$$
$$= 1 + O(n^{1-4\alpha}).$$

It is also standard fare (e.g. [1], Lemma 5.3.3) that

$$\mathbb{E}\{Y\mathbb{1}_{[Y>0]}\} = \frac{\eta}{\sqrt{2\pi}}\exp\left(-\frac{\nu^2}{2\eta^2}\right) + \nu\Phi\left(\frac{\nu}{\eta}\right)$$
$$= \frac{\eta}{\sqrt{2\pi}}\left[1 - \frac{\nu^2}{2\eta^2} + O(n^{2-8\alpha})\right] + \frac{\nu}{2} + \frac{\nu^2}{\sqrt{2\pi}\eta}$$
$$+ O(n^{3/2-6\alpha})$$
(4.13) $$= \frac{1}{\sqrt{2\pi}} + \frac{\nu}{2} + O(n^{\max(-\alpha,3/2-6\alpha)})$$
(4.14) $$= (2\pi)^{-1/2} + O(n^{1/2-2\alpha}).$$

Let $p_k(y) = yr_{k-1}(y)$ and $Y_0 = (Y-\nu)/\eta$. Write $p_k(y) = \sum_{l=1}^k a_{l-1}y^l$, with the $a_i$'s at most $O(n^{-\alpha})$. Then

$$p_k(Y) = \sum_{l=1}^k a_{l-1}(\eta Y_0 + \nu)^l$$
(4.15) $$= \sum_{j=0}^k \left[\sum_{l=\max(1,j)}^k a_{l-1}\binom{l}{j}\eta^j\nu^{l-j}\right]Y_0^j$$
$$= \sum_{j=0}^k c_j Y_0^j,$$

where the largest coefficient $c_j$ in (4.15) can be loosely bounded by $O(n^{-\alpha})$.

The absolute moments of $Y_0$ are all finite, so that

$$|\mathbb{E}p_k(Y)\mathbb{1}_{[Y>0]}| \leq \sum_{j=0}^{k} |c_j|\mathbb{E}|Y_0|^j = O(n^{-\alpha}),$$

whereby

$$\mathbb{E}\{Y(1 + r_{k-1}(Y))\mathbb{1}_{[Y>0]}\} = (2\pi)^{-1/2} + O(n^{1/2-2\alpha}).$$

Putting the pieces together, (4.10) simplifies to

$$e^{-I_u(\hat{\theta}_u)}\tau_u^{-k}(2\pi)^{-(k+1)/2}H_{k-1}(\tau_u\hat{\theta}_u)$$
$$\times (1 + O(n^{1/2-2\alpha}))(1 + O(n^{-\alpha}))(1 + O(n^{1-4\alpha}))$$
$$= \hat{\rho}_{k,n}(u) \times (1 + O(n^{1/2-2\alpha})).$$

Toward (4.6), consider one additional term in (4.13)–(4.14):

$$(2\pi)^{-1/2} + \nu/2 + O(n^{\max\{-\alpha, 3/2-6\alpha\}})$$
$$= (2\pi)^{-1/2} + \dot{\mu}_k/2 + O(n^{\max\{-\alpha, 3/2-6\alpha\}})$$
$$= (2\pi)^{-1/2} + k_3(Z, Z, Z_k)\hat{\theta}_u^2/4 + O(n^{\max\{-\alpha, 3/2-6\alpha\}})$$
$$= (2\pi)^{-1/2}\left\{1 + \mathbb{E}[Z^2 Z_k]\sqrt{\frac{\pi}{2}}\frac{\hat{\theta}_u^2}{2}\right\}(1 + O(n^{\max\{-\alpha, 3/2-6\alpha\}})).$$

As for (4.7), we have

$$-\frac{1}{2}\dot{\mu}|'_{k-1}\dot{\Sigma}|_{k-1}^{-1}\dot{\mu}|_{k-1} = -\frac{1}{2}\dot{\mu}|'_{k-1}(I + O(n^{-\alpha}))\dot{\mu}|_{k-1}$$
$$= -\frac{1}{2}\sum_{i \leq k-1}(\mathbb{E}[Z^2 Z_i]\hat{\theta}_u^2/2 + O(n^{1/2-3\alpha}))^2$$
$$+ \dot{\mu}|'_{k-1}O(n^{-\alpha})\dot{\mu}|_{k-1}$$
$$= -\frac{\hat{\theta}_u^4}{8}\sum_{i \leq k-1}\mathbb{E}[Z^2 Z_i]^2 + O(n^{1-5\alpha}).$$

Finally, note that

$$\max\{1/2 - 3\alpha, 1 - 5\alpha\} < -\alpha \leq \max\{-\alpha, 3/2 - 6\alpha\}$$
$$\leq \max\{-\alpha, 1 - 4\alpha\}. \qquad \square$$

It is worth emphasizing that if the *mixed* (third-order) cumulants $\mathbb{E}[Z^2 Z_i]$ of $Z$ with its first derivatives vanish, then the approximation is improved regardless of the behavior of the *marginal* cumulants of $Z$. Note that if $1/3 \leq$ $\alpha \leq 1/2$, then (4.7) is no better than (4.6), since $-\alpha \geq 1 - 4\alpha$. An important special case of Theorem 4 is when $\alpha = 1/2$, so that $u$ is a constant. In that "small-threshold" setting each of the approximations (4.5)–(4.7) admits a relative error of the order $n^{-1/2}$.

4.2. *What about $\alpha \leq \frac{1}{4}$?* Consider the case $\alpha \leq \frac{1}{4}$, in which the threshold $u$ grows somewhat more quickly relative to $n$. Here we expect the EC heuristic to perform well, and the Gaussian approximation to perform poorly. But what of the tilting procedure? From (3.5) we have that when $\alpha < 1/4$, the tilted mean derivatives also grow with $n$: $\dot\mu_i, \ddot\mu_{ij} \to \infty$. Lemma 3 no longer holds in this scenario since the "remainder" terms now dominate (4.1) and, thus, the special covariance structure which leads to Hermite polynomials is lost (the problem is exacerbated if $\alpha < 0$).

A reasonable conjecture in this setting is that the tilted EC densities will provide a much better approximation to the true ones than will those derived from a Gaussian assumption (see the next section). Less obvious is whether the tilted densities themselves have an acceptable relative error—indeed, Theorem 4 seems to suggest the contrary. While we have yet to work out the details in general, empirical evidence in at least two cases (see Sections 6 and 7) supports using the tilted densities.

4.3. *Brief comment on Rabinowitz and Siegmund.* In the language and notation used here, the Rabinowitz–Siegmund heuristic can be phrased

$$(4.16) \quad \mathbb{P}\Big\{\sup_t Z(t) \geq u\Big\} \approx 1 - \exp\{-\mathbb{E}\{M_u(Z)\}\} \approx \mathbb{E}\{M_u(Z)\},$$

where $M_u(Z)$ is the number of local maxima of the field $Z$ above $u$. In RS, the authors approximated this expectation by

$$(4.17) \quad \mathbb{E}\{M_u(Z)\} \doteq \exp\{-I_u(\hat\theta_u)\} \frac{|T||\widehat\Lambda|}{(2\pi)^{(k+1)/2}\tau_u} \hat\theta_u^{N-1},$$

where $\widehat\Lambda = \mathsf{Var}_{\hat\theta_u}[\dot Z]$ is the tilted variance of the full derivative process (matrix of second spectral moments), having removed the local isotropy condition. Note that $\hat\theta_u^{N-1}/\tau_u$ is the leading term in the polynomial $\tau_u^{-N} H_{N-1}(\tau_u \hat\theta_u)$, so that (modulo the spectral determinant) (4.17) is equivalent to the leading term of our approximation. Naïvely extending (4.17) for, say, a stationary process on a rectangle, one might obtain

$$(4.18) \quad \mathbb{E}\varphi(A_u) \doteq \sum_{k=0}^{N} \sum_{J \in \mathcal{O}_k} |J|_k |\widehat\Lambda_J|^{1/2} \widehat\rho_{k,n}(u),$$

where now $\widehat\Lambda_J$ denotes the tilted variance of the derivatives restricted to facet $J$ of the rectangle, $\mathcal{O}_k$ is the set of $k$-dimensional facets touching the origin (cf. AT, Chapter 11), and $|\cdot|_k$ is a $k$-dimensional Lebesgue measure. (See Section 5.2 for a definition of $\widehat{\rho}_{0,n}$.) For a general isotropic process, where $\widehat{\Lambda} = \widehat{\lambda} I$, (4.18) would become

$$(4.19) \qquad \mathbb{E}\varphi(A_u) \doteq \sum_{k=0}^{N} \widehat{\lambda}^{k/2} \mathcal{V}_k(T) \widehat{\rho}_{k,n}(u).$$

It is readily seen that the relative error between the approximations (4.18)–(4.19) and our own versions (with the untilted $\lambda$ or $\Lambda_J$ replacing $\widehat{\lambda}$ or $\widehat{\Lambda}_J$) is of order $O(n^{-\alpha})$. Indeed, $\lambda_{jl} = \mathbb{E}\{Z_j(t) Z_l(t)\} = O(1)$, and

$$(4.20) \qquad \widehat{\lambda}_{jl} = \sum_{\nu=2}^{\infty} k_\nu(Z, \ldots, Z, Z_j, Z_l) \frac{\hat{\theta}_u^{\nu-2}}{(\nu-2)!} = \lambda_{jl}(1 + O(n^{-\alpha})).$$

Hence, $|\widehat{\Lambda}|^{1/2} = |\Lambda|^{1/2}(1 + O(n^{-\alpha}))$, and similarly for $\widehat{\Lambda}_J$. This implies a $O(n^{-\alpha})$ discrepancy for each $k \geq 1$ term in (4.18) and (4.19). Since $-\alpha \leq 1/2 - 2\alpha$, it is unclear from the theory whether one should use the tilted or untilted second spectral moment(s) in computing $\mathbb{E}\varphi(A_u)$. Similarly, it is seen in Section 5 and the proof of Theorem 4 that using the Gaussian argument $u/\sigma$ in the place of $\tau_u \hat{\theta}_u$ in the Hermite polynomials leaves the relative error unchanged from (4.4). Intuition suggests that $\tau_u \hat{\theta}_u$ may give more accurate results in finite samples. A geometric argument also supports its use; see [11]. The $k=0$ term is unaffected by these changes, and can be improved by the methods to be discussed in Section 5.2.

4.4. *Some remarks about validity.* It has been shown that for a Gaussian process the relative error of the expected EC as an approximation to the $p$-value is exponential as $u \to \infty$, in the sense that the first term of their difference is *exponentially smaller* than the zeroth term in expansion (1.3), which behaves like $u^{-1} e^{-u^2/(2\sigma^2)}$ [33]. The authors of [34] (and an anonymous reviewer) rightly pointed out that this result is, for the time being, a purely Gaussian one: the exponential error likely does not hold even for a $t$ random field. Nevertheless, much of the development in [33] makes no use of Gaussianity. We argue here that for CLRF's, the expected EC is exponentially close to the expected number of local maxima above $u$, which is obviously an upper bound for the $p$-value.

Assuming that the parameter space is a manifold without boundary, one can write the difference (before taking expectations) as

$$(4.21) \quad |\varphi(A_u) - M_u| \leq \#\{t \in T : Z(t) \geq u, \dot{Z}(t) = 0, \ddot{Z}(t) \not< 0\}$$

$$(4.22) \qquad\qquad = \#\{t \in T : Z(t) \geq u, \dot{Z}(t) = 0, (\ddot{Z}(t) + Z(t)I) \not< Z(t)I\}$$

$$(4.23) \qquad\qquad \leq \#\{t \in T : Z(t) \geq u, \dot{Z}(t) = 0, \zeta_1(\ddot{Z}(t) + Z(t)I) > u\},$$

where $\zeta_1(\cdot)$ extracts the largest eigenvalue. We have again assumed $\mathsf{Var}\dot{Z}(t) = I$, purely for convenience. It is not hard to argue (e.g., using AT's Expectation Meta-Theorem in combination with an Edgeworth expansion) that for a stationary CLRF the expectation of (4.23) behaves to first order like

$$|T|e^{-u^2/(2\sigma^2)}e^{-u^2/(2\sigma_c^2)},$$

where $\sigma_c^2 \triangleq \sup_{\|v\|=1} \mathsf{Var}[v'(\ddot{Z}+ZI)v]$ is the so-called critical variance. The key fact is that $\ddot{Z}+ZI$, $\dot{Z}$ and $Z$ are uncorrelated, and indeed *almost* independent. When $T$ has a boundary, (4.21) must account for points where $A_u$ intersects $\partial T$; the argument is more finicky but not fundamentally different. In the interest of brevity, we leave the details to a future article.

Since $p \leq \mathbb{E}M_u \approx \mathbb{E}\varphi(A_u)$, clearly using the expected EC is conservative (at least asymptotically). Determining conditions under which the $p$-value approximation is *sharp* will require modification of the purely Gaussian portions of [33]. Numerical evidence from a simulation (see Section 7) is reassuring.

**5. Comparison with the Gaussian approximation.** We shall use the notation

(5.1) $$\epsilon_n[c_n, d_n] = \beta \quad \Leftrightarrow \quad c_n = d_n \times (1 + O(n^\beta))$$

as $n \to \infty$ to describe relative error for $\beta \leq 0$. We say that the error is *sharp*, and write $\epsilon_n[c_n, d_n] \stackrel{s}{=} \beta$, if (5.1) holds, but fails if $\beta$ is replaced by $\beta - \epsilon$ for any $\epsilon > 0$. It is important to determine whether there are values of $\alpha$ for which our derived expression provides a better approximation to $\rho_{k,n}(u)$ than does the Gaussian EC density. We examine only the approximation $\widehat{\rho}_{k,n}(u)$; similar considerations for the refined formulae (4.6) and (4.7) are straightforward and left to the reader.

5.1. *EC densities with $k \geq 1$.* In Theorem 4 we derived a "tilted" order $O(n^{1/2-2\alpha})$ (in relative error) approximation to $\rho_{k,n}(u)$ for $k \geq 1$, denoted by $\widehat{\rho}_{k,n}(u)$ and given in (4.5). Comparing this approximation to the Gaussian EC density $\rho_k^\gamma = \rho_{k,n}^\gamma$ from (1.6), we have (under the regularity conditions of Theorem 4) the following:

COROLLARY 5. *Let $u \sim cn^{1/2-\alpha}$ for some $0 < c < \infty$. If $\frac{1}{3} \leq \alpha \leq \frac{1}{2}$ and $k_3(Z) \neq 0$, then for all $k \geq 1$,*

(5.2) $$\epsilon_n[\rho_{k,n}^\gamma(u), \rho_{k,n}(u)] \stackrel{s}{=} 1 - 3\alpha.$$

*In particular,*

(5.3) $$\epsilon_n[\rho_{k,n}^\gamma(u), \rho_{k,n}(u)] - \epsilon_n[\widehat{\rho}_{k,n}(u), \rho_{k,n}(u)] \geq \tfrac{1}{2} - \alpha \geq 0.$$

Moreover, if $\frac{1}{4} < \alpha < \frac{1}{3}$ and $k_3(Z) \neq 0$, then the error of the Gaussian densities is exponential, in the sense that

$$\rho_{k,n}^{\gamma}(u) = \rho_{k,n}(u) \times \exp(b_n) \tag{5.4}$$

for $k \geq 1$ and some sequence $b_n$ with $|b_n| \to \infty$.

REMARK. In the proof of Corollary 5, it is easily seen that the sequence referred to in (5.4) can be written $b_n = k_3(Z)\hat{\theta}_u^3/6 + O(n^{1-4\alpha})$, where the skewness may be positive or negative.

5.2. *The $k=0$ case.* Until now we have only considered EC densities for $k \geq 1$. The zeroth density is generally defined by

$$\rho_{0,n}(u) \stackrel{\Delta}{=} \mathbb{P}\{Z(t) \geq u\}, \tag{5.5}$$

the univariate tail probability above the level $u$. Since nothing in the derivations of Theorem 4 or Corollary 5 explicitly requires $u > 0$, the results mentioned thus far are valid if we weaken the threshold condition to $u \sim cn^{1/2-\alpha}$, $|c| < \infty$. [Their use is limited, however, by the facts that (a) when $|u|$ is small the Gaussian approximation is adequate and (b) when $u$ is large and negative the zeroth density dwarfs the others.] With the exception of the Lugananni–Rice formula (5.9), the discussion in this subsection *does* assume $u \geq 0$, since it pertains to saddlepoint-type approximations to (5.5). If one wishes to approximate $\rho_{0,n}(u)$ for $u < 0$ (as is the case in Figure 1), then the expressions can be modified by making use of the symmetry in their derivations and of the fact that $\rho_{0,n}(u) = 1 - \mathbb{P}\{Z(t) < u\}$. Such generalizations are left to the reader. It should be noted that when a $p$-value approximation is sought and hence $u$ is large, the error in $\rho_{0,n}(u)$ typically has a negligible impact compared to those in the higher EC densities. Approximating tail probabilities via saddlepoint methods is far from a new problem; the principal purpose of this subsection is merely to enumerate some useful formulae.

Let $\Psi(\cdot) \stackrel{\Delta}{=} 1 - \Phi(\cdot)$ denote the standard Normal survival function. Given Theorem 4, one who is inspired by AT's definition $H_{-1}(x) = \sqrt{2\pi}\Psi(x)e^{x^2/2}$ might naïvely estimate (5.5) by

$$\hat{\rho}_{0,n}(u) \stackrel{\Delta}{=} \exp\{-I_u(\hat{\theta}_u)\}(2\pi)^{-1/2}H_{-1}(\tau_u\hat{\theta}_u)$$

$$= \exp\{\tau_u^2\hat{\theta}_u^2/2 - I_u(\hat{\theta}_u)\}\Psi(\tau_u\hat{\theta}_u).$$

Perhaps surprisingly, this is what Robinson [26] called the "first" saddlepoint approximation, and turns out to be a reasonable estimate if $u > 0$. The zeroth Gaussian density is

$$\rho_{0,n}^{\gamma}(u) = \Psi\left(\frac{u}{\sigma}\right) = \frac{\sigma e^{-u^2/(2\sigma^2)}}{u\sqrt{2\pi}}(1 + O(n^{2\alpha-1})).$$

Under typical Central Limit Theorem conditions (e.g., [15], Theorem 2.3), one has to first order

$$\rho_{0,n}(u) \doteq \rho_{0,n}^{\gamma}(u) + \frac{k_3(Z)}{6\sigma^3}\left(\frac{u^2}{\sigma^2} - 1\right)\phi(u/\sigma), \quad (5.6)$$

which implies (sharply, assuming nonzero skewness) that

$$\rho_{0,n}(u) = \rho_{0,n}^{\gamma}(u)(1 + O(n^{1-3\alpha})).$$

Daniels [14] gave two differently-derived formulae. The first can be written

$$\rho_{0,n}(u) \doteq \exp\{\tau_u^2 \hat{\theta}_u^2/2 - I_u(\hat{\theta}_u)\}$$
$$\times \left[\Psi(\tau_u\hat{\theta}_u)\left(1 - \frac{K_Z^{(3)}(\hat{\theta}_u)\tau_u^3\hat{\theta}_u^3}{6}\right)\right. \quad (5.7)$$
$$\left. + \phi(\tau_u\hat{\theta}_u)\frac{K_Z^{(3)}(\hat{\theta}_u)}{6}(\tau_u^2\hat{\theta}_u^2 - 1)\right]$$

and corresponds to the so-called "second" saddlepoint approximation [26]. Let $f_Z(\cdot)$ denote the marginal density of $Z(t)$. Upon using the large $x$ estimate $\Psi(x) = \phi(x)/x \times (1 + O(x^{-2}))$, and under the regularity condition (cf. [26]),

$$\sup_z \left|\Phi\left(\frac{z-u}{\tau_u}\right) - \int_{-\infty}^z e^{\hat{\theta}_u y - K_Z(\hat{\theta}_u)} f_Z(y)\,dy\right| = O(n^{-1/2}), \quad (5.8)$$

(5.7) leads to

$$\rho_{0,n}(u) = \hat{\rho}_{0,n}(u)(1 + O(n^{-\alpha})).$$

The second approach is the Lugannani–Rice formula [20]

$$\rho_{0,n}(u) = \Psi(\omega_u) + \phi(\omega_u)\left[\left(\frac{1}{\zeta_u} - \frac{1}{\omega_u}\right) + \frac{1}{\zeta_u}\left(\frac{K_Z^{(4)}(\hat{\theta}_u)}{8} - \frac{5K_Z^{(3)}(\hat{\theta}_u)^2}{24}\right)\right.$$
$$\left. - \frac{K_Z^{(3)}(\hat{\theta}_u)}{2\zeta_u^2} + \left(\frac{1}{\omega_u^3} - \frac{1}{\zeta_u^3}\right) + \cdots\right], \quad (5.9)$$

where $\omega_u = \mathrm{sgn}[\hat{\theta}_u]\sqrt{2I_u(\hat{\theta}_u)}$ and $\zeta_u = \tau_u\hat{\theta}_u$. Note that this approximation is valid for large negative as well as positive $u$, and consequently, the first two terms in (5.9) were used in the astrophysics example of Figure 1. One can reduce (5.9) to

$$\rho_{0,n}(u) = \Psi(\omega_u)(1 + O(n^{-\alpha}))$$
$$= \Psi(\mathrm{sgn}[\hat{\theta}_u]\sqrt{2I_u(\hat{\theta}_u)})(1 + O(n^{-\alpha})) \quad (5.10)$$
$$\stackrel{\Delta}{=} \overline{\rho}_{0,n}(u)(1 + O(n^{-\alpha})).$$

Filling in a few details gives the following:

PROPOSITION 6. *Let* $\frac{1}{4} < \alpha \leq \frac{1}{2}$ *and suppose that (5.8) holds. Then*
$$\epsilon_n[\widehat{\rho}_{0,n}(u), \rho_{0,n}(u)] = \epsilon_n[\overline{\rho}_{0,n}(u), \rho_{0,n}(u)] = -\alpha.$$
*Also, if the field $Z$ is as in (1.1), $k_3(Z) \neq 0$ and $\frac{1}{3} \leq \alpha \leq \frac{1}{2}$, then*
$$\epsilon_n[\rho_{0,n}^\gamma(u), \rho_{0,n}(u)] \stackrel{s}{=} 1 - 3\alpha.$$
*In particular,*
$$\epsilon_n[\rho_{0,n}^\gamma(u), \rho_{0,n}(u)] > \epsilon_n[\widehat{\rho}_{0,n}(u), \rho_{0,n}(u)] \qquad \forall \alpha \in [1/3, 1/2).$$
*When $\alpha < 1/3$ the error in the zeroth Gaussian density diverges to infinity.*

So once again the tilted densities perform better than the Gaussian one. However, $\overline{\rho}_{0,n}$ and $\widehat{\rho}_{0,n}$ can both be improved upon if one is prepared to numerically integrate the pointwise saddlepoint expansion of the density of $z$,

$$(5.11) \quad f(z) = \frac{1}{\sqrt{2\pi}\tau_z} \exp\{-I_z(\hat{\theta}_z)\}\left\{1 + \left[\frac{1}{8}\mu_4(z) - \frac{5}{24}\mu_3(z)^2\right] + \varepsilon_n(z)\right\}.$$

In (5.11) $\varepsilon_n(x)$ is a $O(n^{-2})$ (for fixed $z$) function of $z$, and

$$\mu_3(z) = \tau_z^{-3} K_Z^{(3)}(\hat{\theta}_z); \qquad \mu_4(z) = \tau_z^{-4} K_Z^{(4)}(\hat{\theta}_z)$$

are the standardized skewness and kurtosis of $Z$ under the tilted (conjugate) density. The square-bracketed term in (5.11) is of order $O(n^{-1})$ for each fixed $z$. Hence,

$$(5.12) \quad \begin{aligned} \rho_{0,n}(u) &= \mathbb{P}\{Z(t) \geq u\} \\ &\doteq \int_{z=u}^{\infty} \widehat{f}_n(z)\,dx + \int_{z=u}^{\infty} \widehat{f}_n(z)(\tfrac{1}{8}\mu_4(z) - \tfrac{5}{24}\mu_3(z)^2)\,dz, \end{aligned}$$

where $\widehat{f}_n(z) \stackrel{\Delta}{=} \exp\{-I_z(\hat{\theta}_z)\}/[\sqrt{2\pi}\tau_z]$ is assumed integrable. So provided

$$(5.13) \quad \sup_{z \geq u}|\tfrac{1}{8}\mu_4(z) - \tfrac{5}{24}\mu_3(z)^2 + \varepsilon_n(z)| = o(n^{-\alpha})$$

as $n \to \infty$, the first term in (5.12) will give an approximation with smaller relative error than $\widehat{\rho}_{k,n}(u)$ or $\overline{\rho}_{k,n}(u)$ (assuming numerical integration does not compound the error). The tradeoff is an increase in computation, although solving the saddlepoint equation everywhere can be avoided by a judicious change of variables ([15], Section 6.2). It is well known (e.g., [15], Remark 3.2) that gains can be achieved when $\widehat{f}_n(z)$ is normalized to become a density. Let $\overset{\triangle}{f}_n(z) = \widehat{f}_n(z)/\|\widehat{f}_n\|_1$; a good approximation to $\rho_{0,n}(u)$ is given by

$$\overset{\triangle}{\rho}_{0,n}(u) \stackrel{\Delta}{=} \int_{z=u}^{\infty} \overset{\triangle}{f}_n(z)\,dz.$$

For $p$-value calculation the order $O(n^{-\alpha})$ approximations $\widehat{\rho}_{0,n}$ and $\overline{\rho}_{0,n}$ will often be sufficiently accurate. If the field is not too smooth and/or the search region is large in volume, higher terms begin to dominate (1.3) and the Gaussian density $\rho_{0,n}^\gamma$ may in fact be adequate. For a more complete account of tail probability approximations, see the references already mentioned in this subsection.

Corollary 5 and Proposition 6 suggest that, asymptotically, tilting the EC densities is worthwhile provided $u$ grows with $n$. In statistical applications, "$u \to \infty$" corresponds to "large $u$," which is a crucial if unstated assumption when thresholding random field data via the expected EC heuristic. However, the tilted densities are not guaranteed to be any more accurate than the Gaussian ones in the small threshold scenario $\alpha = 1/2$, at least as $n \to \infty$. The performance of the approximation in finite samples is another matter altogether. These issues and others will be addressed in Sections 6 and 7.

**6. Example: A scaled $\chi_n^2$ random field.** It is instructive to consider a random field for which both the tilted and exact EC densities are available explicitly; the $\chi_n^2$ random field is one such case. Let $X_1, X_2, \ldots, X_n$ be i.i.d. centred, unit variance Gaussian random fields on $T$ and put $Y(t) = \sum_{i=1}^n X_i(t)^2$. It has been established (by two distinct techniques [31, 38]) that the EC densities of $Y$ are given by

$$\rho_{k,n}^\chi(u) = \frac{(k-1)! u^{(n-k)/2} e^{-u/2}}{(2\pi)^{k/2} \Gamma(n/2) 2^{n/2-1}}$$

$$(6.1) \qquad \times \sum_{j=0}^{\lfloor (k-1)/2 \rfloor} \sum_{l=0}^{k-1-2j} \binom{n-1}{k-1-2j-l}$$

$$\times \frac{(-1)^{k-1}(-u)^{j+l}}{l! j! 2^j} \mathbb{1}_{\{n \geq k-2j-l\}}$$

for $k \geq 1$.

We normalize by forming $Z(t) = n^{-1/2}(Y(t) - n)$, so that $Z$ is centered, approximately Gaussian, and $\sigma^2 = \mathsf{Var}[Z(t)] = 2$. Note that the $\chi^2$ field is "twice as rough" as each $X_i$; if the component Gaussian fields are isotropic with second spectral moment $\lambda_x$, then (see [1], Chapter 7) $\mathsf{Var}(\dot{Z}) = n^{-1} \mathsf{Var}(\dot{Y}) = 4\lambda_x I_N$. To insure that $\lambda = 1$ for the normalized process, we therefore assume that the $X_i$'s are isotropic with $\Lambda_X \stackrel{\Delta}{=} \mathsf{Var}(\dot{X}_i) = I_N/4$. With that covariance structure and for a square parameter space $T$, we have

$$\mathbb{E}\varphi(\{t : Z(t) \geq u\}) = \mathbb{E}\varphi(\{t : Y(t) \geq \sqrt{n}u + n\})$$

$$= \sum_k |I_k/4|^{1/2} \rho^\chi_{k,n}(\sqrt{n}u + n)\mathcal{V}_k(T),$$

$$= \sum_k 2^{-k} \rho^\chi_{k,n}(\sqrt{n}u + n)\mathcal{V}_k(T).$$

Thus, the proper comparison is between $\widehat{\rho}_{k,n}(u)$ and $2^{-k}\rho^\chi_{k,n}(\sqrt{n}u + n)$ for $k = 0, \ldots, N$ [this is a minor notational discrepancy; (6.1) is in fact $2^k \rho^\chi_{k,n}(u)$ by the conventions of this article]. We have

$$K_Z(\theta) = -\sqrt{n}\theta - \frac{n}{2}\log\left(1 - 2\frac{\theta}{\sqrt{n}}\right);$$

$$K'_Z(\theta) = \frac{2\theta}{1 - 2\theta/\sqrt{n}};$$

$$K''_Z(\theta) = 2\left(1 - 2\frac{\theta}{\sqrt{n}}\right)^{-2}.$$

It follows that

$$\hat{\theta}_u = \frac{\sqrt{n}u}{2(u + \sqrt{n})}; \qquad I_u(\hat{\theta}_u) = \frac{\sqrt{n}u}{2} - \frac{n}{2}\log\left(1 + \frac{u}{\sqrt{n}}\right);$$

$$\tau_u^2 = 2\left(1 + \frac{u}{\sqrt{n}}\right)^2; \qquad \tau_u\hat{\theta}_u = \frac{u}{\sqrt{2}}.$$

Note that the Hermite portion of the tilted density is identical to that of the Gaussian density in this case. Plugging these into (4.5), we have, for $k \geq 1$,

$$(6.2) \qquad \widehat{\rho}_{k,n}(u) = \frac{e^{-u\sqrt{n}/2}(1 + u/\sqrt{n})^{n/2-k}}{(2\pi)^{(k+1)/2}2^{k/2}}H_{k-1}(u/\sqrt{2}).$$

Collecting all of the constants in (6.1), again for $k \geq 1$,

$$(6.3) \qquad \begin{aligned} \rho_{k,n}(u) &= 2^{-k}\rho^\chi_{k,n}(\sqrt{n}u + n) \\ &= \frac{e^{-u\sqrt{n}/2}(1 + u/\sqrt{n})^{n/2-k}}{(2\pi)^{(k+1)/2}2^{k/2}}R_{k-1,n}(u), \end{aligned}$$

where (assuming $n \geq N$) $R_{k,n}$ is the following degree $k$ polynomial:

$$R_{k,n}(x) = 2^{-(k-2)/2}k!\frac{\sqrt{2\pi}e^{-n/2}}{\Gamma(n/2)}\left(\frac{n}{2}\right)^{(n+1)/2}$$

$$\times \sum_{j=0}^{\lfloor k/2 \rfloor}\sum_{l=0}^{k-2j}\binom{n-1}{k-2j-l}\frac{(-1)^{k+j+l}n^{j+l-k/2-1}}{l!j!2^j}$$

$$\times \left(1 + \frac{x}{\sqrt{n}}\right)^{j+l+(k+1)/2}.$$

Also, note that

$$\frac{\sqrt{2\pi}e^{-n/2}}{\Gamma(n/2)}\left(\frac{n}{2}\right)^{(n+1)/2} \sim \frac{\Gamma(n/2+1)}{\Gamma(n/2)} = \frac{n}{2}$$

as $n \to \infty$. The corresponding Gaussian densities are

(6.4) $$\rho^\gamma_{k,n}(u) = \frac{e^{-u^2/4}}{(2\pi)^{(k+1)/2}2^{k/2}}H_{k-1}(u/\sqrt{2}), \qquad k \geq 1.$$

For the $k=0$ EC density, we have

$$\rho_{0,n}(u) = \rho^\chi_{0,n}(\sqrt{n}u + n) = \mathbb{P}\{Y \geq \sqrt{n}u + n\} = \int_{\sqrt{n}u+n}^{\infty} g_n(y)\,dy,$$

where $Y \stackrel{\mathscr{D}}{=} Y(t) \sim \chi^2_n$ and $g_n(y) = \Gamma(n/2)^{-1}2^{-n/2}y^{n/2-1}e^{-y/2}$ is the density of $Y$. It is a well-known fact [13] that the first-order saddlepoint approximation $\widehat{g}_n$ to $g_n$ is exact up to a constant, so that

(6.5) $$g_n(y) \equiv \stackrel{\triangle}{g}_n(y).$$

To obtain a similar result for $\rho_{0,n}(u)$ and $\stackrel{\triangle}{\rho}_{0,n}(u)$, we shall use the following lemma, which is easily proved.

LEMMA 7. *Let $Z = \sigma_n^{-1}(Y - \mu_n)$ be a random variable with density $f_n(z) = \sigma_n g_n(\sigma_n z + \mu_n)$, where $g_n(y)$ is the density of $Y$. Then $\forall z \in \mathbb{R}$,*

$$\stackrel{\triangle}{f}_n(z) = \sigma_n \stackrel{\triangle}{g}_n(\sigma_n z + \mu_n).$$

Let $f_n(z)$ denote the density of $Z(t)$. Putting $\sigma_n = \sqrt{n}$ and $\mu_n = n$ in the lemma and using (6.5), we have

$$\begin{aligned}\rho^\chi_{0,n}(\sqrt{n}u + n) &= \int_{\sqrt{n}u+n}^{\infty} \stackrel{\triangle}{g}_n(y)\,dy \\ &= n^{-1/2}\int_{\sqrt{n}u+n}^{\infty} \stackrel{\triangle}{f}_n(n^{-1/2}(y-n))\,dy \\ &= \int_u^{\infty} \stackrel{\triangle}{f}_n(z)\,dz \\ &= \stackrel{\triangle}{\rho}_{0,n}(u).\end{aligned}$$

So the integrated, *scaled* tilted zeroth density is exact in this case, and this is a rare instance where normalization of the saddlepoint approximation does not increase computation. The ordinary tilted and Gaussian $k=0$ densities are given, respectively, by

$$\widehat{\rho}_{0,n}(u) = e^{(u^2-2\sqrt{n}u)/4}\left(1 + \frac{u}{\sqrt{n}}\right)^{n/2} \Psi\left(\frac{u}{\sqrt{2}}\right);$$

$$\rho_{0,n}^{\gamma}(u) = \Psi\left(\frac{u}{\sqrt{2}}\right).$$

As a term-by-term comparison of (6.2) and (6.3) is onerous, we compare $H_k(u/\sqrt{2})$ and $R_{k,n}(u)$ for $k=0,1,2$, the cases of interest in a 2D or 3D problem. We have

$$H_0(u/\sqrt{2}) = 1, \quad R_{0,n}(u) = c_n\sqrt{1 + \frac{u}{\sqrt{n}}} = 1 + O(n^{-\alpha});$$

$$H_1(u/\sqrt{2}) = \frac{u}{\sqrt{2}},$$

$$R_{1,n}(u) = c_n \frac{u}{\sqrt{2}}\left(1 + \frac{u}{\sqrt{n}} + \frac{1}{u\sqrt{n}} + \frac{1}{n}\right) = \frac{u}{\sqrt{2}}(1 + O(n^{-\alpha}));$$

$$H_2(u/\sqrt{2}) = \frac{u^2}{2} - 1,$$

$$R_{2,n}(u) = c_n\left(\frac{u^2}{2} - 1 + \frac{u}{2\sqrt{n}} + \frac{1}{n}\right)\left(1 + \frac{u}{\sqrt{n}}\right)^{3/2}$$

$$= \left(\frac{u^2}{2} - 1\right)(1 + O(n^{-\alpha}));$$

where $c_n = \sqrt{2\pi}e^{-n/2}\Gamma(n/2)^{-1}(n/2)^{(n-1)/2} = 1 + O(n^{-1})$. Hence, the tilted EC densities for $k=1,2,3$ have a relative error $O(n^{-\alpha})$. Note that the Gaussianity of the $X_i$'s in this example leads to $\mathbb{E}[Z^2 Z_j] = 0$, which implies that, for all $k \geq 1$, a relative error *no worse* than $O(n^{\max\{-\alpha, 3/2 - 6\alpha\}})$ must be achieved by $\widehat{\rho}_{k,n}$—see (4.7). It is easily verified that $I_u(\widehat{\theta}_u) = u^2/4 - u^3/(6\sqrt{n}) + O(n^{1-4\alpha})$, so that $\epsilon_n[\rho_{k,n}^{\gamma}(u), \rho_{k,n}(u)] \stackrel{s}{=} 1 - 3\alpha$ when $\alpha \geq 1/3$.

The behavior of these approximations for $k \leq 3$ is explored further in Figure 2. In each column we consider a different asymptotic regime: the most general one in which $\alpha = 1/3$ and both $n$ and $u$ grow; the classical one in which $\alpha = 1/2$ so that $u$ is fixed while $n$ grows; and the less common one in which $n$ is fixed while $u$ grows (corresponding to $\alpha = -\infty$). The latter regime is meaningful for statistical applications since the sample size $(n)$ is large but finite and the relative error in $\mathbb{E}\varphi(A_u)$ as a $p$-value approximation is conjectured to fall off exponentially as $u \to \infty$ (see Section 4.4). The theory of Sections 4.1 and 5 does not directly apply to this scenario. A fourth regime, where $k \to \infty$, might also be of interest; empirical evidence for the $\chi_n^2$ field suggests that the tilted densities perform poorly under that scheme.

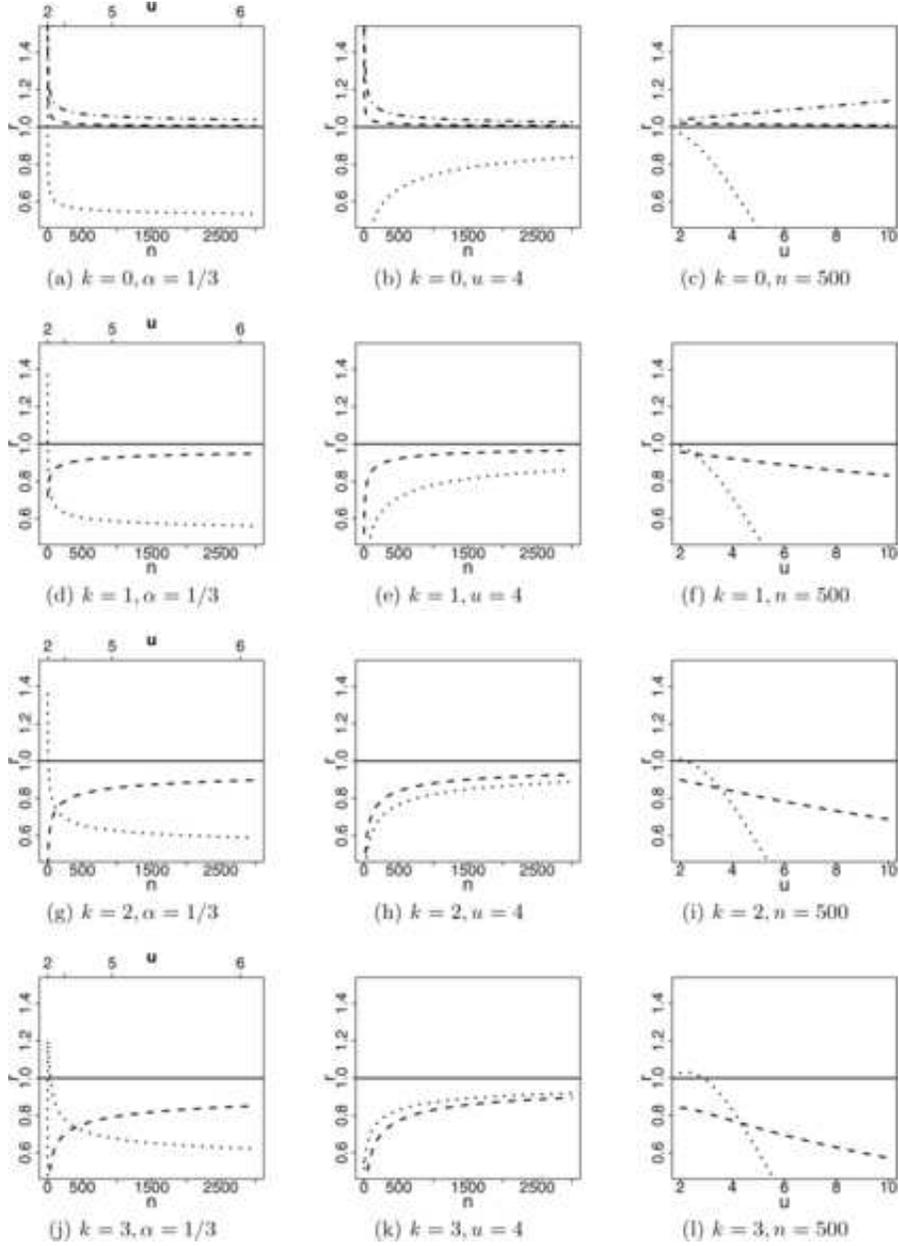

FIG. 2. *EC densities for the $\chi_n^2$ random field: exact ($\rho_{k,n}$, "—"); tilted ($\widehat{\rho}_{k,n}$, "- -"); and Gaussian ($\rho_{k,n}^\gamma$, "···"). Ratios $r$ of approximate to exact densities are plotted for three regimes: $u \sim cn^{1/6} \to \infty$ ($\alpha = 1/3$, 1st col.); $u$ fixed and $n \to \infty$ ($\alpha = 1/2$, 2nd col.); $n$ fixed and $u \to \infty$ ($\alpha = -\infty$, 3rd col.). For $k = 0$, $\rho_{0,n} \equiv \overset{\triangle}{\rho}_{0,n}$, and "-··-" denotes $\overline{\rho}_{0,n}$.*

For $k = 0, 1$ in this example, the tilted densities offer a clear improvement over the Gaussian approximation even for small $n$ and/or $u$. For $k = 2, 3$, a larger sample size or threshold is needed before the difference emerges, but it is still unequivocal when $u \to \infty$. When $\alpha = 1/2$ for $k = 2$ or $3$, the tilted and Gaussian approximations perform similarly, as predicted by Corollary 5. In the third regime where $u \to \infty$ while $n$ remains fixed, the Gaussian densities are quickly diverging from the correct densities; however, they are more accurate than the tilted densities for $k = 2, 3$ in the range $2 \leq u \leq 4$. Note that the chosen sample size $n = 500$ was quite large in that regime. When $k = 0$ the naïve tilted density $\widehat{\rho}_{0,n}(u)$ appears to outperform the Lugannani–Rice estimate's leading term $\overline{\rho}_{0,n}(u)$, but this difference is small and may be specific to the $\chi_n^2$ field. The Gaussian assumption generally leads to *under* estimation of the expected EC in this example. In other words, using a Gaussian assumption for statistical thresholding of $\chi_n^2$ fields is anti-conservative.

**7. Application: Subtracted "bubbles" images.** A moment's thought reveals that the form (1.1) taken for the $Z$ process is far more restrictive than is necessary. Indeed, the crucial property of $\mathbf{Z} = (Z, \dot{Z}', \ddot{Z}')$ appearing in the proof of Theorem 4 is the rate of decay of its (joint) cumulants; namely

(7.1) $$\mathbb{E}[Z(t)] \equiv 0, \qquad \mathsf{Var}[Z(t)] \equiv \sigma^2;$$

and

(7.2) $$k_\nu([\mathbf{Z}(t)]^\gamma) = O(n^{-(\nu-2)/2}), \qquad \gamma \in \mathbb{N}^{(N+1)(N+2)/2}, |\gamma| = \nu, \nu \geq 2.$$

One might therefore expect Theorem 4 to apply to a variety of stationary random fields satisfying (7.1) and (7.2) in addition to the regularity conditions in the Appendix. The following is a slightly more general construction than (1.1): let $\{W_i(t)\}_i$ be independent but perhaps nonidentically-distributed random fields with finite cumulants of all orders. Assume that the random field of interest takes the form

(7.3) $$Z(t) = \frac{1}{\sqrt{n}} \sum_{i=1}^n a_i W_i(t)$$

for finite constants $(a_i)_{i=1}^\infty$. This framework covers the case of a two-sample problem (e.g., a difference in means between two groups), as well as the following application.

The bubbles methodology of Gosselin and Schyns has been described in detail elsewhere [16]. Typically, the experimenters' goal is to determine which parts of an image are most important in a visual discrimination task. The data analyzed later in this section come from [3]. On each trial of the experiment, a subject was shown a $256^2$ image of a face that was either fearful or happy. The images were masked apart from a random number of localized regions or "bubbles" that revealed only selected parts of the face; see Figure 3. The subject was then asked whether the partially revealed face was fearful or happy, and the trial was repeated about 3,000 times on each of 10 subjects.

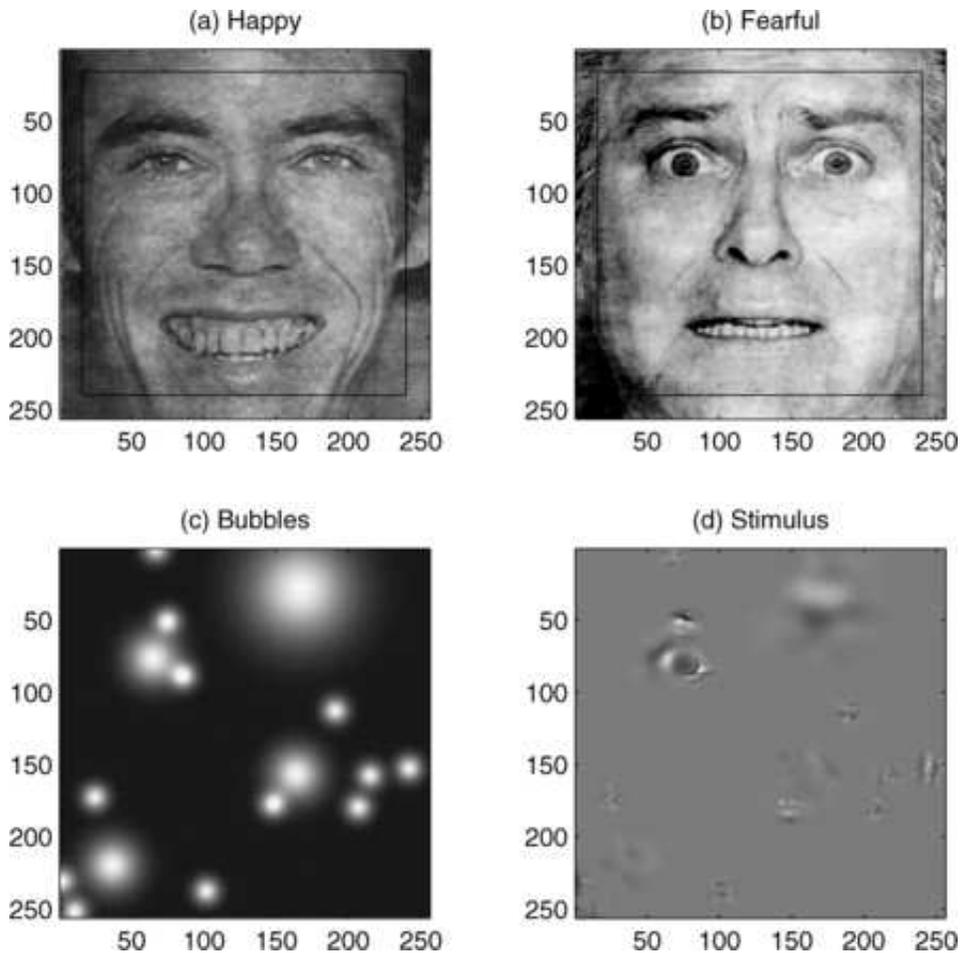

Fig. 3. *Bubbles experiment. The subject is asked to discriminate between the happy* (a) *and fearful* (b) *faces on presentation of the stimulus* (d) *which is one of the two faces, here the fearful face, partially revealed by random "bubbles"* (c). *The $227^2$ search region $T$ is inside the black frame in* (a) *and* (b).

Each masked image was generated as follows. First a range of scales was chosen, $\sigma_{Bj} = 3 \times 2^j$, $j = 1, \ldots, 5$. Then the original image $\mathcal{I}_0$ [Figure 3(a) or (b)] was smoothed by an isotropic Gaussian filter with standard deviation $\sigma_{Bj}/2$ to produce images $\mathcal{I}_j$. The smoothed images were differenced to produce images $D_j = \mathcal{I}_{j-1} - \mathcal{I}_j$ that reveal image features at five different scales. Differenced image $D_j$ was then multiplied by a mask consisting of the sum of a random number of isotropic Gaussian "bubbles," each with standard deviation $\sigma_{Bj}$. The bubble center were chosen at random from the $256^2$ pixels, that is, according to a Poisson process. The number of bubbles for each scale was a multinomial random variable with probabilities inversely proportional to bubble area ($\propto \sigma_{Bj}^{-2}$), so that equal areas were revealed by each of the five bubble sizes, on average. The sum of all the bubbles over all five scales is shown in Figure 3(c) for one of the trials. The sum of all the bubbles times the differenced images is shown in Figure 3(d). On the basis of Figure 3(d) the subject must decide if the face is happy or fearful (fearful is the correct answer in this case). The total number of bubbles was chosen dynamically to maintain a roughly 0.75 probability of correctly identifying the face.

7.1. *Test statistic, cumulants and tilted densities.* For the present analysis we consider a single-subject run. We assume that the grid is square, and that the bubble standard deviation $\sigma_B$ (or equivalently, the bubble full width at half max $F = \sqrt{8 \log 2}\sigma_B$) is fixed, as is the number $m$ of bubbles per image. We denote by $p_C$ and $p_I = 1 - p_C$ the proportions of correctly and incorrectly classified images. Let $Y_{ip}$ denote the number of bubbles in image $i$ which are centred at pixel $p$, with $i = 1, \ldots, n$ and $p = 1, \ldots, P$. Then $Y_{ip} \sim \mathcal{P}(m/P)$ is a Poisson random variable, and the $Y_{ip}$ are approximately independent (exactly when they come from different images). The test statistic field for detecting significant regions is formed by signing and summing the bubbles, with sign depending on whether the image was correctly or incorrectly classified. Thus, the unnormalized test statistic at pixel $t$ is given by

$$(7.4) \qquad \widetilde{Z}(t) \triangleq n^{-1/2} \sum_{i=1}^{n} \sum_{p=1}^{P} a_i b_p(t) Y_{ip},$$

where

$$a_i = \begin{cases} 1/p_C, & \text{image } i \text{ correct,} \\ -1/p_I, & \text{image } i \text{ incorrect,} \end{cases}$$

and

$$b_p(t) = \exp\left\{-\frac{\|p - t\|^2}{2\sigma_B^2}\right\}.$$

The statistic is normalized by forming

$$Z(t) = \widetilde{Z}(t)/\sqrt{\mathsf{Var}\widetilde{Z}(t)}.$$

Inference can be either conditional on the observed vector $a$ or unconditional, in which case the $a_i$ are i.i.d. scaled Rademacher random variables and $p_C$ is taken as the (assumed known) *population* proportion. One may apply the theory of the present article under either analysis; we have chosen the conditional approach since it simplifies cumulant calculations. Thus, $p_C$ and $p_I$ are observed proportions.

The most natural global null hypothesis is

$$H_0 : a_i \text{ is independent of } (Y_{ip})_{p=1}^{P} \qquad \forall i \leq n,$$

which implies

$$H_0' : \mathbb{E}\{Z(t)|a\} = 0 \qquad \forall t \in T,$$

since $\mathbb{E}_{H_0}\{Z(t)|a\} \propto \sum_p b_p(t) \sum_i a_i \mathbb{E}_{H_0}\{Y_{ip}|a\} = \frac{m}{P}(\sum_p b_p(t))(\sum_i a_i) = 0$. Ignoring mild boundary effects, note that the null distribution of $Z(t)$ is unchanged if we translate or rotate the coordinate system; thus, $Z(t)$ is (strictly) isotropic. By (7.4), $Z$ is of the form $n^{-1/2}\sum_i a_i W_i$ where the $W_i$'s are i.i.d. random fields. The distribution of $W_i(t) = \sum_p b_p(t) Y_{ip}$ is very non-Gaussian, with a high probability of being nearly zero (if there are no bubble center near $t$), and a small probability of being large. The distributions become more Gaussian as $m$ or $F$ increase.

An important practical note is that most users of this methodology choose to "clamp" the summed bubbles at a maximum height equal to the peak height of a single bubble. In particular, two bubbles centered at the same pixel, when added, have the same peak height as a single bubble. This is how the mask is actually created, and it is felt that this is how the stimulus is perceived. Clamping effectively replaces (7.4) with $n^{-1/2}\sum_{i=1}^{n} a_i \min\{\sum_{p=1}^{P} b_p(t)Y_{ip}, 1\}$. Doing so not only complicates cumulant derivations, but, more importantly, degrades the "smoothness" of the random field. When the number of bubbles is small compared to the number of pixels and their width is not too large (as in the data analyzed below), the loss of smoothness due to clamping is negligible. Our choice of an unclamped analysis is always valid, but perhaps less powerful at detecting facial features if indeed the probability of response is related to the clamped rather than the unclamped bubble mask.

Without any extra effort we generalize the cumulant computations to $N$ dimensions, although currently in bubbles experiments $N = 2$. It can be shown that the $j$th null conditional cumulant of $Z(t)$ is approximately

$$(7.5) \qquad \kappa_j = \left(\frac{2^{j/2}}{j}\right)^{N/2} \frac{p_C^{1-j} + (-1)^j p_I^{1-j}}{(p_C^{-1} + p_I^{-1})^{j/2}} \left\{\left(\frac{4\log 2}{\pi}\right)^{N/2} \frac{1}{B}\right\}^{j/2-1},$$

which only depends on $p_C$ and the total number $B \triangleq nm/(PF^{-N})$ of bubbles per $N$-dimensional "resel" (volume element). These have been used to compute the tilted variance of $Z$,

$$\tau_u^2 = \sum_{j=0}^{\infty} \kappa_{j+2} \frac{\hat{\theta}_u^j}{j!}.$$

In practice (for typical values of $p_C$ and $B$), $j \leq 20$ seems to be enough; beyond this the cumulants are extremely close to zero. The tilting parameter $\hat{\theta}_u$ is found by solving $K'_Z(\theta) = u$ with a Newton–Raphson algorithm started from the Gaussian solution $\theta = u$; five iterations appear to suffice in most cases. An approach like this one was advocated by RS.

One can also show that the second spectral moment of $Z$ under $H_0$ is given (approximately) by

(7.6) $$\mathsf{Var}\{\dot{Z}(t)|a\} \approx \frac{1}{2\sigma_B^2} I_N \triangleq \lambda I_N.$$

The $p$-value approximation when $N = 2$ is thus

(7.7) $$\mathbb{P}\left\{\sup_t Z(t) \geq u\right\} \approx \mathbb{E}\varphi(A_u) \doteq \hat{\rho}_0(u) + 2(P\lambda)^{1/2}\hat{\rho}_1(u) + P\lambda\hat{\rho}_2(u),$$

where the EC densities are

(7.8) $$\begin{aligned}\hat{\rho}_0(u) &= \exp\{\tau_u^2\hat{\theta}_u^2/2 - I_u(\hat{\theta}_u)\}\Psi(\tau_u\hat{\theta}_u),\\ \hat{\rho}_1(u) &= \exp\{-I_u(\hat{\theta}_u)\}(2\pi\tau_u)^{-1},\\ \hat{\rho}_2(u) &= \exp\{-I_u(\hat{\theta}_u)\}(2\pi)^{-3/2}\tau_u^{-1}\hat{\theta}_u.\end{aligned}$$

By contrast, the Gaussian theory EC densities are

(7.9) $$\begin{aligned}\rho_0^\gamma(u) &= \Psi(u),\\ \rho_1^\gamma(u) &= \exp\{-u^2/2\}(2\pi)^{-1},\\ \rho_2^\gamma(u) &= \exp\{-u^2/2\}(2\pi)^{-3/2}u,\end{aligned}$$

with $p$-value approximation

(7.10) $$\rho_0^\gamma(u) + 2(P\lambda)^{1/2}\rho_1^\gamma(u) + P\lambda\rho_2^\gamma(u).$$

When $P$ is large, the third terms in (7.7) and (7.10) dominate. Typical values for the parameters are $P = 256^2$, $m = 16.5$, $F = 14.1$, $n = 3000$, $p_C = 0.75$; setting $u = 3.965$ [the $p = 0.05$ threshold using (7.7)] then gives

$$\hat{\rho}_0(u) = 1.4 \times 10^{-5}, \qquad 2(P\lambda)^{1/2}\hat{\rho}_1(u) = 0.00171, \qquad P\lambda\hat{\rho}_2(u) = 0.0484.$$

In view of this, the ratio of the tilted to Gaussian $p$-value approximations,

(7.11) $$r \approx \frac{\hat{\rho}_2(u)}{\rho_2^\gamma(u)} = \exp\left\{\frac{u^2}{2} - I_u(\hat{\theta}_u)\right\}\frac{\hat{\theta}_u}{\tau_u u},$$

only depends on $u$, $B$ and $p_C$. We also consider the approximation suggested by extension (4.18) of RS, namely,

$$\widehat{\rho}_0(u) + 2(P\widehat{\lambda}_u)^{1/2}\widehat{\rho}_1(u) + P\widehat{\lambda}_u\widehat{\rho}_2(u), \tag{7.12}$$

where $\widehat{\lambda}_u$ denotes the *tilted* second spectral moment of $Z$. It can be demonstrated that $\widehat{\lambda}_u \approx \tau_u^2 \lambda$. To compare this approximation to the Gaussian theory, one can again form the $p$-value ratio

$$r_{\text{RS}} \approx \frac{\widehat{\lambda}_u \widehat{\rho}_2(u)}{\lambda \rho_2^{\gamma}(u)} \approx \tau_u^2 r. \tag{7.13}$$

These ratios are plotted for $u = 3.5$ and a range of $B$ and $p_C$ in Figure 4. Note that in this application the approximation (7.7) deviates *less* extremely from the Gaussian one than does the heuristic of RS. This can be understood by noting that, when $\alpha = 1/2$,

$$\tau_u^2 = 1 + \kappa_3 \hat{\theta}_u + O(n^{-1}), \tag{7.14}$$

where $\hat{\theta}_u \approx u > 0$ and $\text{sgn}\{\kappa_3\} = \text{sgn}\{\frac{1}{2} - p_C\}$. At the same time, $r > 1 \Leftrightarrow p_C \leq 1/2$. The $O(n^{-1})$ term in (7.14) can be positive or negative but is of smaller order than the skewness term. It is also worth noting that unlike in the $\chi_n^2$ example and the fields considered by RS, the Gaussian approximation can be conservative (i.e., lead to larger $p$-values) in bubbles experiments when typical parameter values ($p_C > 1/2$) are used.

7.2. *Simulation.* A simulation was devised to compare (7.7), (7.10) and (7.12) with the true tail probability (1.2) of $\sup_t Z(t)$ above a threshold under $H_0$. In each of $10^5$ Monte Carlo iterations, we generated 3000 bubbles images on a $P = 208 \times 208$ square with corresponding classifications, using $p_C = 0.75$ and bubble width $F = 24$. For each iteration we then calculated a sequence of test statistic fields $Z_n(t)$ for $n = 100, 200, \ldots, 3000$. Since kernel smoothing with Fourier methods was used to speed up computation, the data were generated first on a larger $256 \times 256$ square before discarding the outer edges to eliminate periodicity. The number of bubbles per image was fixed so that an average of $m = 20$ would fall within the search region; this has a negligible effect on the cumulants (7.5). The mixed skewness terms $\mathbb{E}[Z(t)^2 Z_i(t)]$ can be shown to be (approximately) proportional to $\sum_p (t_i - p_i) b_p(t)^3$, and hence, are very close to zero for pixels far from the boundary (recall that they were exactly zero in the $\chi^2$ example). We therefore expect the tilted approximations to be within the order $O(n^{-\alpha})$ of the true $p$-value. Figure 5 displays the simulation results under the three asymptotic regimes. In addition to using both $\lambda$ and $\widehat{\lambda}_u$, we also tried $u/\sigma$ in the place of $\tau_u \hat{\theta}_u$ as the argument in the Hermite polynomials. This resulted in third and fourth variations on the tilted $p$-value approximation. In all instances the tilted expected EC's were closer to the observed $p$-value than the Gaussian expected EC. In particular, the expansion which uses the tilted spectral moment $\widehat{\lambda}_u$ and the Gaussian Hermite argument $u/\sigma$ ("– –," orange) performed best, although differences between tilted approximations were small. The observed directional trend between the four tilted expressions is likely specific to this application. We note that while feasible for parameter spaces of this size, deriving $p$-values by simulation for bubbles-type data is time-consuming, and quickly becomes impracticable as the dimension $N$ increases beyond 2. The benefits of having a closed-form $p$-value approximation are evident.

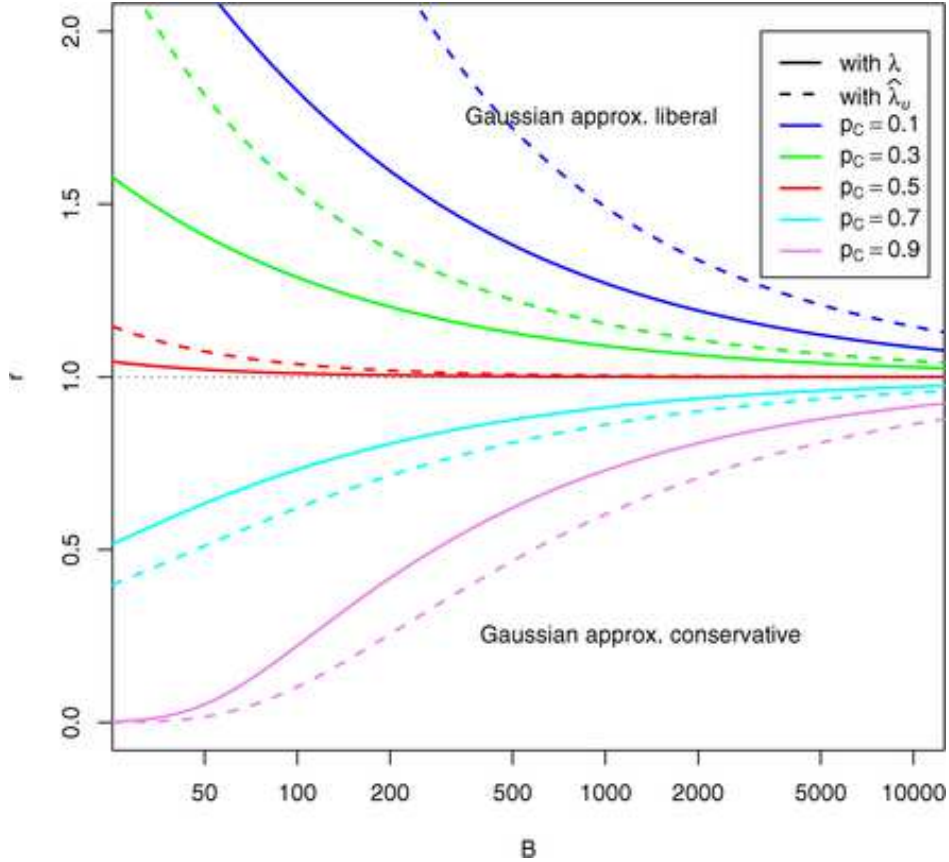

Fig. 4. *Ratios of tilted to Gaussian-approximated expected Euler Characteristic of the excursion of a test statistic from a typical bubbles experiment above the threshold $u = 3.5$, using both the untilted ($\lambda$) and the tilted ($\widehat{\lambda}_u$) second spectral moment. The ratios are plotted for a reasonable range of $B$ (the total number of bubbles per resel) and for several values of $p_C$ (the proportion of correctly classified images). When $p_C = 1/2$ (red curve), the skewness of $Z$ is zero, and the saddlepoint correction is less pronounced. Assuming that the tilted expressions are closer to the true p-value, $r > 1$ signifies that the Gaussian approximation is anti-conservative, while $r < 1$ means that the Gaussian procedure is conservative. Note that $p_C \approx 0.75$ is commonly targeted by the originators of the bubbles paradigm, making the Gaussian theory valid though less powerful than the tilted inference.*

7.3. *Data analysis.* We analyzed a portion of the data from [3] described at the beginning of this section (control group). One subject (#9) was presented with $n = 3072$ images of size $256 \times 256$ pixels, which he/she classified as either happy or fearful, at a rate of $p_C = 0.76$. The search region was trimmed to the inner $P = 227^2$ pixels to avoid Fourier edge effects. In agreement with the theory in this article, we have limited our analysis to a single bubble width, $F = 18.8$ pixels or $\sigma_B = 8$ (the subject was in fact shown bubbles of width $F = 28.3$, i.e., $j = 2$, $\sigma_{Bj} = 12$; but we chose a narrower filter when smoothing the test statistic in order to totally eliminate the "clamping" issue). The mean number of bubbles per image at this scale and within the search region was $m = 3.0$. Our chosen single bubble width was the second-narrowest class ($j = 2$); an example of the three bubbles is clearly visible in Figure 3(c). These parameter values correspond to $B = 63.6$ bubbles/resel. Thresholded images are displayed in Figure 6. There is strong signal in the right eye region for this subject. The global maximum of $Z(t)$ exceeded 5.8, and hence, all thresholding methods identified that peak. At the $p = 0.05$ level the Gaussian threshold picked up only the right eye, while the tilted threshold found additional activation in the upper lip. At $p = 0.1$ significance the tilted expected EC identified a third signal in the left eye which went undetected by the Gaussian theory. Note that the Gaussian approximation is conservative in this case since $Z(t)$ exhibits negative skewness ($p_C > 1/2$).

**8. Discussion and future work.** We have provided a rigorous theoretical justification for the use of saddlepoint methods in thresholding asymptotically-Gaussian random fields, first proposed by Rabinowitz and Siegmund [25]. In particular, we have extended their "expected number of maxima" heuristic to an expected Euler Characteristic approximation for locally isotropic Central Limit-type fields, and derived relative errors for the EC densities when the threshold is allowed to grow. We have also compared these formulae to a Gaussian approximation, both analytically and in numerical and real data examples involving a $\chi_n^2$ field and a linear combination of Poisson fields.

As did RS (who only considered a Poisson process) vis-à-vis expected local maxima, we found substantial differences between the Gaussian and the various tilted approximations (which were uniformly more accurate in a particular case—see Figure 2) to the expected EC. It was shown in Section 5 that in a large-threshold setting the saddlepoint correction reduces asymptotic relative error. In contrast to the examples chosen by RS, we have demonstrated via simulation that the Gaussian approximation can be conservative in some instances rather than always underestimating the $p$-value; see Section 8 and Figures 4–5. We have shown that the important correction for non-Gaussianity occurs in the exponential portion of the EC densities. Additional corrections to the variance of the first derivative field and to the argument of the Hermite polynomial are of smaller order and may not be needed unless the expected EC over the full range of $u$ is required. The examples considered indicate that tilting the spectral moment has a greater impact than tilting the Hermite argument. Our results support RS's recommendation to employ a tilted expansion when the additional computation is feasible; we add that the benefit is greatest if the field has nonzero skewness and/or if $u \to \infty$ (is large) along with $n$.

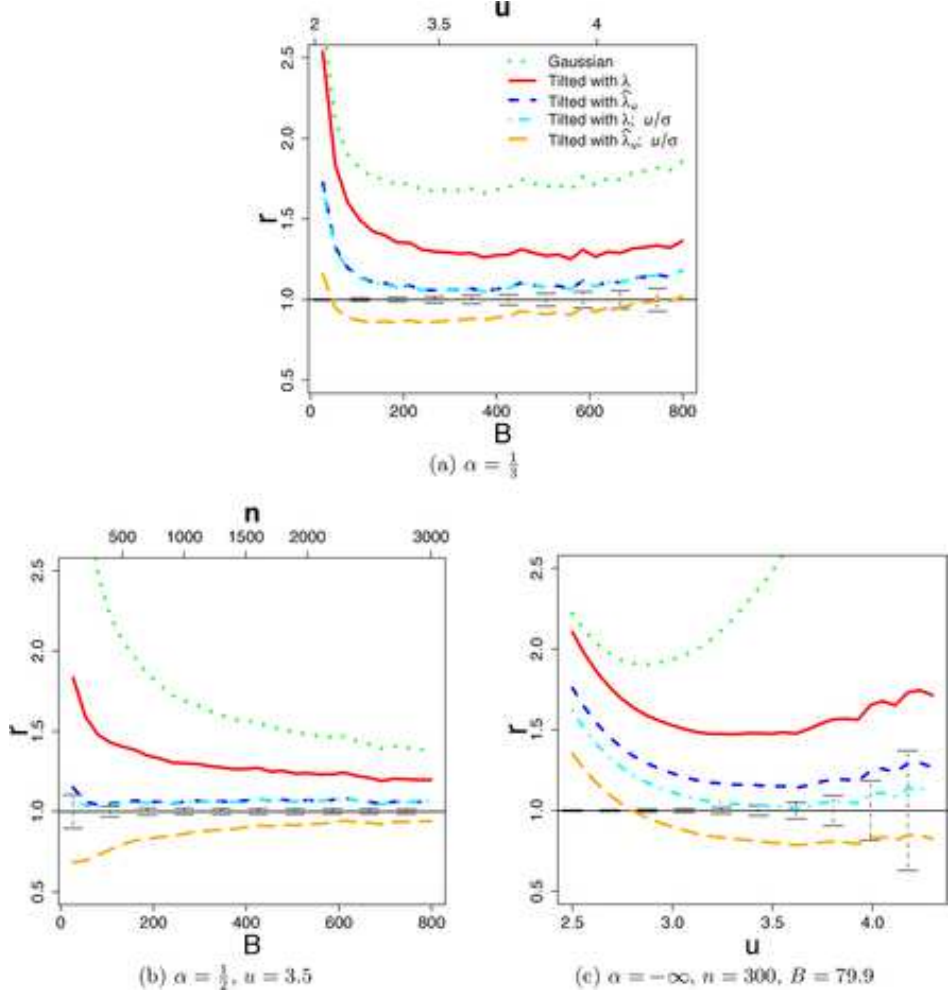

FIG. 5. *Ratios of expected EC approximations to true (empirical) p-values for simulated bubbles data, under the three asymptotic regimes of Section 7. Bubbles parameter values were $P = 208^2, F = 24, m = 20, p_C = 0.75$, and $n$ ranging from 100 to 3000. The ratio scale was chosen so that absolute deviations from the line $r = 1$ correspond to relative error. The Gaussian expected EC is given by (7.10). The tilted approximations use either the untilted ($\lambda$) or tilted ($\widehat{\lambda}_u$) second spectral moment, and either the tilted ($\tau_u \hat{\theta}_u$) or untilted ($u/\sigma$) Hermite argument. Monte Carlo error bars ($\pm 1$ standard deviation) are for the Gaussian curve ("$\cdots$"), but are attached to the line $r = 1$ to aid interpretation. They are computed pointwise, using $\mathsf{Var}(\widehat{p}/p_{emp}) \approx 10^{-5} \widehat{p}^2 (1 - p_{emp})/p_{emp}^3$, where $p_{emp}$ is an empirical p-value and $\widehat{p}$ one of the EC approximations. Error bars for the other curves are slightly narrower. The threshold $u = 3.5$ was chosen in* (b) *in order to have empirical probabilities near 0.05. The sample size $n = 300$ was chosen in* (c) *in order to have the total bubbles/resel close to the value $B = 63.6$ observed in the data. In all cases the tilted approximations outperform the Gaussian one, including in* (c) *where $n$ is fixed, $u$ grows and all curves appear to diverge from unity.*

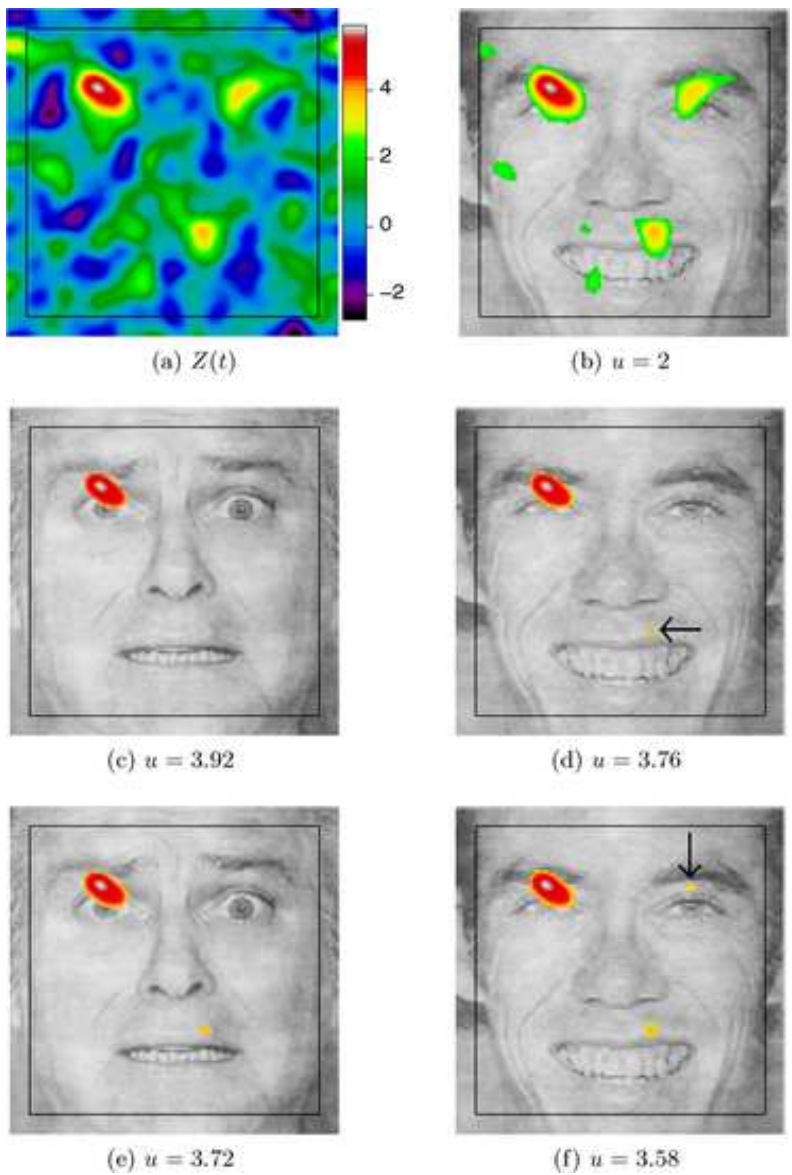

Fig. 6. *Thresholded bubbles data for subject #9, with the search region framed. Test statistic* (a), *with excursions above:* (b) $u = 2$ *(Gaussian uncorrected $p \approx 0.05$);* (c) $u = 3.92$ *(Gaussian $p \approx 0.05$);* (d) $u = 3.76$ *(tilted $p \approx 0.05$);* (e) $u = 3.72$ *(Gaussian $p \approx 0.1$);* (f) $u = 3.58$ *(tilted $p \approx 0.1$). P-values approximated by $\mathbb{E}EC$ under the Gaussian (7.10) and tilted theory with argument $\tau_u \hat{\theta}_u$ and spectral moment $\widehat{\lambda}_u$; see (7.12). Inference using (7.7) is similar and omitted. There is strong signal in the right eye region. At $p = 0.05$, the Gaussian theory identifies only the right eye* (c); *tilting finds activation on the upper lip* (d). *At $p = 0.1$, tilting detects signal above the left eye* (f); *the Gaussian approximation does not* (e). *Tilting reduces thresholds because $p_C > 1/2$, so $Z(t)$ exhibits negative skewness.*

In addition to the bubbles test statistic field presented above, we have begun to apply the theory of this article to lesion density maps from neuroimaging [10], which can be non-Gaussian in a similar fashion. Handling this and further applications will require generalizing these results to certain nonisotropic and nonstationary CLRF's; this work is already underway. As discussed in Section 4.4, it would also be desirable to rigorously extend the proof in [33] of the EC heuristic to the asymptotically Gaussian case. In a forthcoming paper we shall explore a geometric interpretation to the tilted EC densities in the spirit of [31]. This will tie into similar investigations of "CL-related" processes, including asymptotically $\chi^2$, $t$ and $F$ random fields.

## APPENDIX: REGULARITY CONDITIONS FOR THEOREM 4

We assume that the random vectors $\mathbf{W}^i = (W^i, \dot{W}^i|_k', \text{vec}[\sigma \ddot{W}^i|_{k-1} + W^i I/\sigma]')'$ are nonlattice. (While saddlepoint expansions have been established for lattice distributions, the resulting random fields could scarcely be imagined to be regular in the manner described below.) There are simple conditions which validate the tilted expansion (2.6) of the joint density $f(z, \dot{z}, m)$ of $(Z(t), \dot{Z}|_k(t)', M(t)')$ [which exists, and is bounded and continuous, under said conditions]. For example, it is sufficient that: (i) $\mathbf{W}^i$ come from an exponential family indexed by $\theta \in \Theta$ an open, convex subset of $\mathbb{R}^{d_k}$ [which appears tacitly in (2.2)]; and (ii) $\mathbf{W}^i$ have characteristic function $\psi_\theta$ for which $|\psi_\theta(s)|^\nu$ is integrable for some $\nu = \nu(\theta) \geq 1$. For the latter, compare [4], Theorem 6.5, and their condition $(c')$. An extension of Theorem 4 to the case where the $W_i$'s are nonidentically distributed could also be conceived; for regularity conditions the reader is referred to Theorem 6.6 of the same monograph.

The exponential family assumption guarantees that the $\mathbf{W}^i$ have finite joint cumulants of all orders, but to bound the error we shall assume, in particular, that for every $\nu \leq 4$ the cgf of $\mathbf{W}^i$ belongs to $C^\nu(V)$ for some neighborhood $V$ containing 0. We shall also require that for some finite constants $H > 0, (b_i)_{i \geq 2}$ and all $s \geq 2$,

$$\text{(A.1)} \qquad \mathbb{E}\|\mathbf{W}^i\|^s \leq \frac{s!}{2} b_i^2 H^{s-2}.$$

As for the validity of (1.5), conditions are laid out in [1] and AT as follows:

A $T \in \mathbb{R}^N$ is compact and $\partial T$ has zero Lebesgue measure.
B $Z$ is almost surely *suitably regular* on $T$ at the level $u$, namely:
  1. $Z$ has continuous partial derivatives up to second order in an open neighborhood of $T$;
  2. there is no $t \in T$ such that $Z(t) - u = Z_1(t) = \cdots = Z_k(t) = 0$;
  3. there is no $t \in \partial T$ and permutation $\rho$ of $\{1, \ldots, k\}$ such that
  $$Z(t) - u = Z_{\rho_1}(t) = \cdots = Z_{\rho_{k-1}}(t) = 0;$$
  4. there is no $t \in T$ and permutation $\rho$ of $\{1, \ldots, k\}$ such that
  $$Z(t) - u = Z_{\rho_1}(t) = \cdots = Z_{\rho_{k-1}}(t) = \det(Z_{\rho_i \rho_j}(t))_{i,j \leq k-1} = 0.$$

C The second-order derivatives of $Z$ have finite variances.
D The joint density of $(Z, \dot{Z}|'_k, \ddot{Z}|'_{k-1})$ is continuous in each argument.
E The conditional density of $(Z, \dot{Z}|'_{k-1})$ given $(Z_k, \ddot{Z}|'_{k-1})$ is bounded above.
F The moduli of continuity $\omega_i$ and $\omega_{ij}$ of the various derivatives satisfy

$$\mathbb{P}\left\{\max_{i,j \leq k}[\omega_i(\eta), \omega_{ij}(\eta)] > \epsilon\right\} = o(\eta^k)$$

as $\eta \to 0$ for every $\epsilon > 0$.

Note that under the oft-made Gaussian assumption, many of these conditions are trivially satisfied.

**Acknowledgments.** The authors thank two anonymous reviewers. A portion of this research was completed while Nicholas Chamandy was a guest of the Department of Statistics, Stanford University.

N. Chamandy
Department of Mathematics
 and Statistics
McGill University
805 Sherbrooke O.
Montréal, Quebec
Canada H3A 2K6
E-mail: chamandy@math.mcgill.ca
and
Google
1600 Amphitheatre Pkwy.
Mountain View, California 94043
USA
E-mail: chamandy@google.com

K. J. Worsley
Department of Mathematics
 and Statistics
McGill University
805 Sherbrooke O.
Montréal, Quebec
Canada H3A 2K6
E-mail: keith.worsley@mcgill.ca
and
Department of Statistics
University of Chicago
5734 S. University Anevue
Chicago, Illinois 60637
USA

J. Taylor  
Department of Statistics  
Stanford University  
Stanford, California 94305  
USA  
E-mail: johathan.taylor@stanford.edu

F. Gosselin  
Department of Psychology  
Université de Montréal  
Montréal, Quebec  
Canada H3C 377  
E-mail: frederic.gosselin@umontreal.ca